\documentclass[journal,transmag]{IEEEtran}

\usepackage{graphicx}
\usepackage{xcolor}
\usepackage{hyperref}
\usepackage{amsmath,amssymb}
\usepackage{tikz,circuitikz}
\usetikzlibrary{calc}
\usepackage{pgfplots}
\usepackage{tuda-pgfplots}
\usepackage{pgfplotstable}
\pgfplotsset{compat=1.18}
\usepackage[exponent-product = \cdot]{siunitx}
\usepackage{subcaption}

\usetikzlibrary{calc,decorations.markings,shapes.misc,positioning,matrix,shapes.geometric}
\usetikzlibrary{patterns,patterns.meta}
\usetikzlibrary {quotes, arrows.meta}
\usetikzlibrary {graphs}
\usetikzlibrary {graphs.standard}
\usetikzlibrary{arrows.meta,matrix}
\usetikzlibrary{positioning}

\newcommand{\CURL}{\ensuremath{\operatorname{curl}}}
\newcommand{\DIV}{\ensuremath{\operatorname{div}}}
\newcommand{\GRAD}{\ensuremath{\operatorname{grad}}}
\newcommand{\permittivity}{\ensuremath{\epsilon}}
\newcommand{\reluctivity}{\ensuremath{\nu}}
\newcommand{\conductivity}{\ensuremath{\sigma}}
\newcommand{\complexConductivity}{\ensuremath{\kappa}}
\newcommand{\imNum}{\ensuremath{i}}
\newcommand{\Afield}{\ensuremath{\boldsymbol{A}}} %
\newcommand{\Bfield}{\ensuremath{\boldsymbol{B}}} %
\newcommand{\Dfield}{\ensuremath{\boldsymbol{D}}} %
\newcommand{\Efield}{\ensuremath{\boldsymbol{E}}} %
\newcommand{\Jfield}{\ensuremath{\boldsymbol{J}}} %
\newcommand{\rhofield}{\ensuremath{\varrho}} %
\newcommand{\phifield}{\ensuremath{\varphi}} %

\newcommand{\nfield}{\ensuremath{\boldsymbol{n}}}
\newcommand{\zerofield}{\ensuremath{\boldsymbol{0}}}

\newcommand{\diffV}{\ensuremath{\mathrm{d}\mathrm{V}}}
\newcommand{\complConj}[1]{\ensuremath{\overline{#1}}}
\newcommand{\scalarWeakSpace}{\ensuremath{H^1(\Omega)}}%
\newcommand{\vectorWeakSpace}{\ensuremath{H(\Omega;\CURL)}}%

\newcommand{\trans}{\ensuremath{^{\!\top}}}
\newcommand{\scalarBaseFun}{\ensuremath{v}}
\newcommand{\vectorBaseFun}{\ensuremath{\boldsymbol{w}}}
\newcommand{\stiffMatCurl}[1]{\ensuremath{\mathbf{K}_{#1}}}
\newcommand{\stiffMatGrad}[1]{\ensuremath{\mathbf{K}}_{#1}}
\newcommand{\massMatCurl}[1]{\ensuremath{\mathbf{M}_{#1}}}
\newcommand{\gradMat}[1]{\ensuremath{\mathbf{G}_{#1}}}
\newcommand{\divMat}[1]{\ensuremath{\mathbf{S}_{#1}}}
\newcommand{\scalarDOF}{\ensuremath{u}}
\newcommand{\scalarDOFs}{\ensuremath{\mathbf{\scalarDOF}}}
\newcommand{\vectorDOF}{\ensuremath{a}}
\newcommand{\vectorDOFs}{\ensuremath{\mathbf{\vectorDOF}}}
\newcommand{\norm}[2]{\left\lVert#1\right\rVert_{#2}}
\newcommand{\scalarDiscreteSpace}{\ensuremath{V}}%
\newcommand{\vectorDiscreteSpace}{\ensuremath{W}}%

\makeatletter
\renewenvironment{cases}[1][l]{\matrix@check\cases\env@cases{#1}}{\endarray\right.}
\def\env@cases#1{%
  \let\@ifnextchar\new@ifnextchar
  \left\lbrace\def\arraystretch{1.2}%
  \array{@{}#1@{\quad}l@{}}}
\makeatother

\tikzset{
	cuboid/.pic={
		\tikzset{%
			every edge quotes/.append style={midway, auto},
			/cuboid/.cd,
			#1
		}
		\draw [every edge/.append style={pic actions, densely dashed, opacity=.5}, pic actions]
		(0,0,0) coordinate (o) -- ++(-\cubescale*\cubex,0,0) coordinate (a) -- ++(0,-\cubescale*\cubey,0) coordinate (b) edge coordinate [pos=1] (g) ++(0,0,-\cubescale*\cubez)  -- ++(\cubescale*\cubex,0,0) coordinate (c) -- cycle
		(o) -- ++(0,0,-\cubescale*\cubez) coordinate (d) -- ++(0,-\cubescale*\cubey,0) coordinate (e) edge (g) -- (c) -- cycle
		(o) -- (a) -- ++(0,0,-\cubescale*\cubez) coordinate (f) edge (g) -- (d) -- cycle;
	},
	/cuboid/.search also={/tikz},
	/cuboid/.cd,
	width/.store in=\cubex,
	height/.store in=\cubey,
	depth/.store in=\cubez,
	units/.store in=\cubeunits,
	scale/.store in=\cubescale,
	width=10,
	height=10,
	depth=10,
	units=cm,
	scale=.1,
}

\tikzset{
	annotated cuboid/.pic={
		\tikzset{%
			every edge quotes/.append style={midway, auto},
			/cuboid/.cd,
			#1
		}
		\draw [every edge/.append style={pic actions, densely dashed, opacity=.5}, pic actions]
		(0,0,0) coordinate (o) -- ++(-\cubescale*\cubex,0,0) coordinate (a) -- ++(0,-\cubescale*\cubey,0) coordinate (b) edge coordinate [pos=1] (g) ++(0,0,-\cubescale*\cubez)  -- ++(\cubescale*\cubex,0,0) coordinate (c) -- cycle
		(o) -- ++(0,0,-\cubescale*\cubez) coordinate (d) -- ++(0,-\cubescale*\cubey,0) coordinate (e) edge (g) -- (c) -- cycle
		(o) -- (a) -- ++(0,0,-\cubescale*\cubez) coordinate (f) edge (g) -- (d) -- cycle;
		\path [every edge/.append style={pic actions, |-|}]
		(b) +(0,-5pt) coordinate (b1) edge ["22 \cubeunits"'] (b1 -| c)
		(b) +(-5pt,0) coordinate (b2) edge ["22 \cubeunits"] (b2 |- a)
		(c) +(3.5pt,-3.5pt) coordinate (c2) edge ["22 \cubeunits"'] ([xshift=3.5pt,yshift=-3.5pt]e)
		;
	},
	/cuboid/.search also={/tikz},
	/cuboid/.cd,
	width/.store in=\cubex,
	height/.store in=\cubey,
	depth/.store in=\cubez,
	units/.store in=\cubeunits,
	scale/.store in=\cubescale,
	width=10,
	height=10,
	depth=10,
	units=cm,
	scale=.1,
}

\usepackage{eso-pic}
\AddToShipoutPicture*{
	\footnotesize\sffamily\raisebox{0.8cm}{\hspace{1.4cm}\fbox{
		\parbox{\textwidth}{
			© 2025 IEEE. Personal use of this material is permitted. Permission from IEEE must be obtained for all other uses, including reprinting/republishing this material for advertising or promotional purposes, collecting new collected works for resale or redistribution to servers or lists, or reuse of any copyrighted component of this work in other works.
		}
	}}
}

\begin{document}

\title{A Two-Step Formulation of Maxwell's Equations
Using Generalized Tree-Cotree Gauges for Low-Frequency-Stability}

\author{\IEEEauthorblockN{Leon Herles\IEEEauthorrefmark{1},
Mario Mally\IEEEauthorrefmark{1,2},
Jörg Ostrowski\IEEEauthorrefmark{3}, and
Sebastian Schöps\IEEEauthorrefmark{1},
Melina Merkel\IEEEauthorrefmark{1}}
\IEEEauthorblockA{\IEEEauthorrefmark{1}Computational Electromagnetics Group, Technische Universität Darmstadt, 64289 Darmstadt, Germany}
\IEEEauthorblockA{\IEEEauthorrefmark{2}Department of Applied Mathematics, Universidade de Santiago de Compostela, 15782, Santiago de Compostela, Spain}
\IEEEauthorblockA{\IEEEauthorrefmark{3}Siemens Digital Industries Software, Switzerland}%
\thanks{Received 19 February 2025; revised 23 April 2025; accepted 15 May 2025.
Date of publication 29 May 2025.
Corresponding author: Leon Herles (email: leon.herles@stud.tu-darmstadt.de). %
}}

\markboth{Generalized Tree-Cotree Gauges for Low-Frequency-Stability}%
{Generalized Tree-Cotree Gauges for Low-Frequency-Stability}

\IEEEtitleabstractindextext{%
\begin{abstract}
This paper presents a new low-frequency stabilization for a two-step formulation solving the full set of Maxwell's equations. The formulation is based on a electric scalar and magnetic vector potential equation using the electroquasistatic problem as gauge condition. The proposed stabilization technique consists of an adequate frequency scaling for the electroquasistatic problem, and a tree-cotree decomposition of the magnetic vector potential such that its divergence remains consistent with the partial decoupling of the magnetic and electric potentials. The paper discusses two variants and demonstrates the effectiveness by a few computational examples.
\end{abstract}

}

\maketitle

\IEEEdisplaynontitleabstractindextext

\IEEEpeerreviewmaketitle

\section{Introduction}
\IEEEPARstart{M}{any} applications require accurate simulation results of electromagnetic devices over a wide frequency range, e.g., to access electromagnetic compatibility. This commonly includes static, low frequency and high frequency regimes. In such cases, static or quasi-static approximations are insufficient to capture the full functionality of the device.
Therefore, we focus here on formulations that consistently solve the complete set of Maxwell's equations, allowing the evaluation of a device in both the low and high frequency regimes.

A popular starting point is the electric field formulation, e.g. \cite[Section 1.2.1]{Monk_2003aa}. Its finite element discretization is straightforward and leads to symmetric linear systems of equations. However, the condition number of such systems tends to infinity as the frequency tends to zero, see \cite{Zhu_2010aa}. Various approaches have been proposed, all of which decompose the field, typically based on a tree-cotree decomposition \cite{Albanese_1990aa,Munteanu_2002aa}, and its governing equations to maintain control in the frequency limit. They either start directly from the electric field or decompose it into the magnetic vector potential and the electric scalar potential first, see e.g. \cite{Dyczij-Edlinger_1999aa,Hiptmair_2008aa,Jochum_2015aa,Jochum_2016aa,Eller_2017aa, Ho_2016aa}. The main differences are in the treatment of lossy materials, the symmetry and the conditioning of the linear system as well as the number of degrees of freedom, e.g., how many Lagrange multipliers are introduced.

In this contribution we investigate the two-step formulation of Maxwell's equations originally proposed in \cite{Ostrowski_2021aa}. It decouples fields and equations on the continuous level such that existing solvers can be reused when considering adequate source terms. However, we reconsider its low-frequency stability and propose a novel tree-cotree decomposition strategy.

The paper is structured as follows. In \autoref{sec:two-step} the two-step formulation is recapitulated and discretized in \autoref{sec:discrete}. \autoref{sec:lf-stab} discusses low-frequency stabilization techniques and their application to the two-step method. Then, \autoref{sec:tests} shows the results of a numerical reference implementation. The paper closes with conclusions in \autoref{sec:conclusions}.

\section{Formulation}\label{sec:two-step}
Let us consider the two-step approach from \cite{Ostrowski_2021aa} in frequency domain, i.e.,
\begin{align}
    -\DIV\left(\complexConductivity\GRAD\phifield\right) &=\imNum\omega\rhofield^{s} \label{eq:strong1}\\
    \CURL\left(\reluctivity\CURL\Afield\right) + \imNum\omega\complexConductivity\Afield  &= \Jfield^{\mathrm{s}}-\complexConductivity\GRAD\phifield \label{eq:strong2}
\end{align}
in its strong form on a bounded, open and simply-connected domain $\Omega\subset\mathbb{R}^3$. In \eqref{eq:strong1}, one solves an electroquasistatic (EQS) problem with given charge density $\rhofield^{s}$, electric scalar potential $\phifield$ and the complex conductivity $\complexConductivity=\conductivity+\imNum\omega\permittivity$ which depends on angular frequency $\omega=2\pi f$, with $f$ being the ordinary frequency, conductivity $\conductivity\geq0$ and permittivity $\permittivity>0$. The solution to \eqref{eq:strong1} is then used to compute the right-hand-side of \eqref{eq:strong2} which consists of an electricquasistatic source term $-\complexConductivity\GRAD\phifield$ and eventually a source current density,
\begin{equation}
    \Jfield^{\mathrm{s}}=\Jfield^{\mathrm{s}}_0+\imNum\omega\Dfield^{\mathrm{s}}_{\mathrm{e}}~~\text{s.t.}~\DIV\Jfield^{\mathrm{s}}_0=0~~\text{and}~\DIV\Dfield^{\mathrm{s}}_{\mathrm{e}}=-\rhofield^{s}.\label{eq:divProperties}
\end{equation}
The magnetic vector potential $\Afield$ solves \eqref{eq:strong2} where $\reluctivity>0$ denotes the reluctivity. Typical fields of interest are the magnetic flux density $\Bfield=\CURL\Afield$, the electric displacement field
\begin{equation*}
    \Dfield=\permittivity\Efield=\underbrace{\left(-\permittivity\GRAD\phifield\right)}_{=\Dfield_{\mathrm{e}}}+\underbrace{\left(-\imNum\omega\permittivity\Afield\right)}_{=\Dfield_{\mathrm{m}}} + \Dfield^{\mathrm{s}}_{\mathrm{e}}
\end{equation*}
and the %
current density
\begin{equation*}
    \Jfield=\conductivity\Efield=\underbrace{\left(-\conductivity\GRAD\phifield\right)}_{=\Jfield_{\mathrm{e}}}+\underbrace{\left(-\imNum\omega\conductivity\Afield\right)}_{=\Jfield_{\mathrm{m}}} + \Jfield^{\mathrm{s}},
\end{equation*}
where electroquasistatic contributions are indicated by subscript $\mathrm{e}$, and full-Maxwell
corrections by subscript $\mathrm{m}$.
The conductivity $\conductivity$ vanishes in the nonconducting (air) region denoted by $\Omega_{\mathrm{A}}\subset\Omega$; its support is denoted by $\Omega_{\mathrm{C}}\subset\Omega$.
Let $\Gamma$ be the interface between both regions.
Furthermore, we require $\operatorname{supp}(\Jfield^{\mathrm{s}})\subset\Omega_{\mathrm{A}}$ as well as $\operatorname{supp}(\rhofield^{\mathrm{s}})\subset\Omega_{\mathrm{A}}$ for the respective excitations.

Appropriate boundary conditions need to be provided to obtain a well-posed problem. In this work, we use combinations of homogeneous Dirichlet ($\phifield=0$; $\Afield\times\nfield=\zerofield$) and homogeneous Neumann ($\partial\phifield/\partial\nfield=0$; $\CURL\Afield\times\nfield=\zerofield$) boundary conditions (BCs), where $\nfield$ is the normal vector on the boundary. For imposing inhomogeneous Dirichlet BCs the usual approaches can be used \cite[A2]{Formaggia_2012aa}.

\subsection{Inductive Coupling}\label{sec:Darwin}
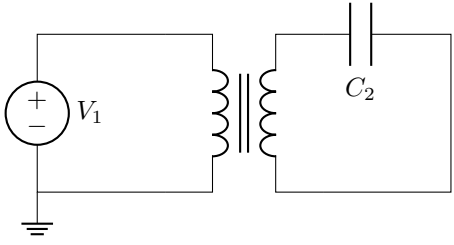
\begin{figure}
    \centering
    \begin{circuitikz}[american]
        \draw (0,0) node [transformer core](T){}
              (T.A1)
              (T.A2)
              (T.B1)
              (T.B2)
              (T.base) node{};
        \draw (T.A1) to ($(T.A1)-(1,0)$) -- ($(T.A1)-(1.7,0)$) to[V, l=$V_1$] ($(T.A2)-(1.7,0)$) -- (T.A2);
        \draw (T.B2) -- ($(T.B2)+(1.7,0)$) -- ($(T.B1)+(1.7,0)$) --($(T.B1)+(1,0)$) to[C, l=$C_2$] (T.B1);
        \draw ($(T.A2)-(1.7,0)$) node[ground]{};
    \end{circuitikz}
    \caption{Example configuration where capacitive and inductive effects are separate but must be considered simultaneously}
    \label{fig:circuit}
\end{figure}
The two-step formulation is a reformulation of the full set of Maxwell's equations and therefore includes all electromagnetic effects. %
Nonetheless, the decomposition into potentials
may lead to counter-intuitive situations. Let us consider the field model of a series connection of a  transformer and a capacitor, visualized as a circuit in \autoref{fig:circuit}. Let us further assume that the model is excited by $V_1=\sin(2\pi f t)\si{\volt}$. Then, the electroquasistatic solution of \eqref{eq:strong1} will yield a current flow and the left leg of the transformer but not through the right one since the inductive coupling is not yet considered. Consequently, the electric field and displacement currents due to $C_2$ will only be visible after solving \eqref{eq:strong2}. This requires the correction term $\Dfield_{\mathrm{m}}$ of the displacement field, which would not be present in Darwin formulations~\cite{Clemens_2022aa}.

This example demonstrates that one shall not attribute capacitive effects to one of the two equations or potentials alone (while this works for inductive effects). Consequently, one must always solve both equations, even if one is possibly interested in capacitive behavior only.

\subsection{Low-Frequency Breakdown}
It is clear that neither \eqref{eq:strong1} nor \eqref{eq:strong2} are uniquely solvable in the static limit $\omega=0$. The electroquasistatic problem \eqref{eq:strong1} degenerates into a stationary current flow problem but includes non-conductive domains \cite{Balian_2023aa}. On the other hand, \eqref{eq:strong2} suffers from the usual non-uniqueness of the curl operators also known from the electric field formulation \cite{Zhu_2010aa} and Darwin formulation~\cite{Kaimori_2024aa}, i.e., gradient fields can be arbitrarily added to the solution. This causes the low-frequency breakdown when numerically solving the problem.

\section{Discretization}\label{sec:discrete}
Following \cite{Ostrowski_2021aa}, we employ the weak formulation: Find $\phifield\in\scalarWeakSpace$ and $\Afield\in\vectorWeakSpace$ such that
\begin{equation}
	\int_{\Omega} \left(\complexConductivity\GRAD\phifield\right)\cdot\GRAD\complConj{\phifield'}\diffV = \imNum\omega\int_{\Omega}\rhofield^{s}\complConj{\phifield'}\diffV \label{eq:weak1}
\end{equation}
for all $\phifield'\in\scalarWeakSpace$ and
\begin{equation}
    \begin{aligned}
        \int_{\Omega} \left(\reluctivity\CURL\Afield\right)\cdot\CURL\complConj{\Afield'} + \imNum\omega\complexConductivity\Afield\cdot\complConj{\Afield'}\diffV \\ = \int_{\Omega}\left(\Jfield^{\mathrm{s}}-\complexConductivity\GRAD\phifield\right)\cdot\complConj{\Afield'}\diffV \label{eq:weak2}
    \end{aligned}
\end{equation}
for all $\Afield'\in\vectorWeakSpace$. The spaces $\scalarWeakSpace$ and $\vectorWeakSpace$ must be endowed with the respective BCs for which we assume that they are chosen such that \eqref{eq:weak1}-\eqref{eq:weak2} is well-posed.

\subsection{Matrices and Block Structures}
\label{sec:discreteMatrices}
We follow the usual Galerkin approach. Testing and approximating equations \eqref{eq:weak1}-\eqref{eq:weak2} and potentials with the same basis functions from finite (sub-)spaces, i.e., $\scalarBaseFun\in\scalarDiscreteSpace\subset\scalarWeakSpace$ and $\vectorBaseFun\in\vectorDiscreteSpace\subset\vectorWeakSpace$. This leads to the matrices
\begin{align}
    \left(\stiffMatGrad{\star}\right)_{ij} &= \int_{\Omega}\star(\GRAD \scalarBaseFun_j)\cdot(\GRAD \scalarBaseFun_i)\,\diffV, \\
    \left(\gradMat{\star}\right)_{ij} &= \int_{\Omega}\star(\GRAD \scalarBaseFun_j)\cdot\vectorBaseFun_i\,\diffV, \label{eq:gradMat} \\
    \left(\massMatCurl{\star}\right)_{ij} &= \int_{\Omega}\star\vectorBaseFun_j\cdot\vectorBaseFun_i\,\diffV,\\
\intertext{where $\star\in\{\conductivity,\permittivity,\complexConductivity\}$ denotes the employed material and the curl-curl matrix}
\left(\stiffMatCurl{\reluctivity}\right)_{ij} &= \int_{\Omega}\reluctivity
    (\CURL \vectorBaseFun_j)\cdot(\CURL \vectorBaseFun_i)\,\diffV.
\end{align}
Together they form the linear equation system
\begin{align}
    \stiffMatGrad{\complexConductivity}\scalarDOFs &= \imNum\omega\mathbf{q}_{\mathrm{s}} \label{eq:discrete1}\\
    \mathbf{W}\vectorDOFs &= \mathbf{j}(\mathbf{u})\label{eq:discrete2}
\end{align}
with $\mathbf{W}=\stiffMatCurl{\reluctivity}+\imNum\omega\massMatCurl{\complexConductivity}$
and $\mathbf{j}(\scalarDOFs) = \mathbf{j}_{\mathrm{s}} - \gradMat{\complexConductivity}\scalarDOFs$. The degrees of freedom (DOFs) of the scalar and vector potential are given as $\scalarDOFs\in\mathbb{C}^{n_{\mathrm{v}}}$ and $\vectorDOFs\in\mathbb{C}^{n_{\mathrm{w}}}$, respectively. The source contributions $\mathbf{q}_{\mathrm{s}}$ and $\mathbf{j}_{\mathrm{s}}$ are constructed as
\begin{align*}
    \left(\mathbf{q}_{\mathrm{s}}\right)_{i} &= \int_{\Omega} \rhofield^{s} \scalarBaseFun_i\,\diffV, \\
    \left(\mathbf{j}_{\mathrm{s}}\right)_{i} &= \int_{\Omega} \Jfield^{s} \cdot \vectorBaseFun_i\,\diffV.
\end{align*}
The discrete spaces $\scalarDiscreteSpace$ and $\vectorDiscreteSpace$ have to satisfy the usual properties of nodal and edge element spaces. An interested reader is referred to \cite{Monk_2003aa}. Here, we chose the spaces based on multipatch isogeometric analysis (IGA) as introduced in \cite{Buffa_2011aa}.
Let us define the two index sets, cf. \cite{Balian_2023aa}
\begin{align}
\mathcal{I}^{(\mathrm{A})}_{\mathrm{v}} &= \{1 \leq n \leq n_{\mathrm{v}} \mid \operatorname{supp}(v_n) \cap\Omega_\mathrm{C}=\emptyset, \, v_n \in \scalarDiscreteSpace\}
\\
\mathcal{I}^{(\mathrm{C})}_{\mathrm{v}} &= \{1 \leq n \leq n_{\mathrm{v}} \mid n \notin \mathcal{I}^{(\mathrm{A})}_{\mathrm{v}}\}
\end{align}
which decompose the finite set of nodal basis functions into the ones supported only in the nonconducting domain and the others (with at least partial support in the conductor). We define $\mathcal{I}^{(\mathrm{A})}_{\mathrm{w}}$ and $\mathcal{I}^{(\mathrm{C})}_{\mathrm{w}}$ analogously for the edge functions. This gives rise to a block matrix partitioning, in particular
\begin{align}
    \massMatCurl{\complexConductivity}
    =
    \begin{bmatrix}
        \massMatCurl{\conductivity}^{(\mathrm{CC})} & \textbf{0} \\
        \textbf{0} & \textbf{0} \\
    \end{bmatrix}
    +
    i\omega
    \begin{bmatrix}
        \massMatCurl{\permittivity}^{(\mathrm{CC})} & \massMatCurl{\permittivity}^{(\mathrm{CA})} \\
       \massMatCurl{\permittivity}^{(\mathrm{AC})} &
       \massMatCurl{\permittivity}^{(\mathrm{AA})} \\
    \end{bmatrix}
\end{align}
where $\massMatCurl{\conductivity}^{(\mathrm{CC})}$ is (symmetric) positive definite while all other blocks of $\massMatCurl{\conductivity}$ vanish due to $\conductivity=0$ in $\Omega_{\mathrm{A}}$.

\subsection{Low-Frequency Breakdown on Matrix Level}
For small $\omega$, i.e., for $\omega\rightarrow 0$, the strong formulation \eqref{eq:strong1}-\eqref{eq:strong2} and consequently its discretization \eqref{eq:discrete1}-\eqref{eq:discrete2} suffer from a low-frequency breakdown. For \eqref{eq:discrete1}, this becomes obvious when examining
\begin{equation}
    \stiffMatGrad{\complexConductivity} = \begin{bmatrix}
        \stiffMatGrad{\conductivity}^{(\mathrm{CC})} + \imNum\omega\stiffMatGrad{\permittivity}^{(\mathrm{CC})} & \imNum\omega \stiffMatGrad{\permittivity}^{(\mathrm{CA})} \\
        \imNum\omega\stiffMatGrad{\permittivity}^{(\mathrm{AC})} & \imNum\omega\stiffMatGrad{\permittivity}^{(\mathrm{AA})} \\
    \end{bmatrix},
\end{equation}
which is singular for $\omega=0$. This also implies that the conditioning deteriorates with decreasing $\omega$. The ratio $\frac{\conductivity}{\permittivity}\gg1$ in $\Omega_\mathrm{C}$ increases this issue even further. Effective modifications are proposed in \cite{Balian_2023aa} to improve the condition number. %
The focus of our paper is the stabilization of \eqref{eq:discrete2} which can be written as
\begin{equation}
    \mathbf{W}=\stiffMatCurl{\reluctivity}-\begin{bmatrix}
       \omega^2\massMatCurl{\permittivity}^{(\mathrm{CC})} - \imNum\omega\massMatCurl{\conductivity}^{(\mathrm{CC})} &  \omega^2\massMatCurl{\permittivity}^{(\mathrm{CA})} \\
        \omega^2\massMatCurl{\permittivity}^{(\mathrm{AC})} &  \omega^2\massMatCurl{\permittivity}^{(\mathrm{AA})} \\
    \end{bmatrix}.
\end{equation}
Here, the contributions from $\imNum\massMatCurl{\complexConductivity}$ vanish for $\omega=0$ and only the singular $\stiffMatCurl{\reluctivity}$ remains. Its kernel consists of discrete gradient fields, and additional issues arise because of large scaling differences due to $\omega\conductivity\gg\omega^2\permittivity$ for $\omega\rightarrow0$.

\section{Low-Frequency Stabilization Techniques}\label{sec:lf-stab}
In order to combat this low-frequency breakdown, we need to recover the generalized Coulomb gauge that has been implicitly enforced by considering the two-step approach \eqref{eq:strong1} and \eqref{eq:strong2}. It becomes evident in the strong form when a divergence is applied to \eqref{eq:strong2} and divide by $\imNum\omega$, yielding the implicit constraint
\begin{equation}
    \DIV\left(\complexConductivity\Afield\right) = 0, \label{eq:strongGenCoul}
\end{equation}
where we exploited that $\DIV\left(\Jfield^{\mathrm{s}}-\complexConductivity\GRAD\phifield\right)=0$ due to \eqref{eq:strong1} and \eqref{eq:divProperties}. This gauge apparently vanishes inside non-conducting regions if $\omega=0$ and numerically already deteriorates in the discrete setting %
if $\omega$ approaches zero \cite{Ostrowski_2021aa}.

This motivates us to split \eqref{eq:strongGenCoul} into
\begin{equation}
    \begin{cases}[r]
        \beta\DIV\left(\left(\conductivity + \imNum\omega\permittivity\right)\Afield\right)=0 & \text{in}~~\Omega_{\mathrm{C}} \\
        \gamma\DIV\left(\permittivity\Afield\right)=0 & \text{in}~~\Omega_{\mathrm{A}}
    \end{cases}\label{eq:modGenCoul}
\end{equation}
which has frequency-independent components in both regions.
The scaling factors $\beta,\gamma\in\mathbb{R}_{>0}$ are added to further reduce the condition number but are assumed to be $\beta=\gamma=1$ if not stated otherwise. The discrete representation of \eqref{eq:modGenCoul} is given as
\begin{equation}
    \divMat{}\vectorDOFs=\begin{bmatrix}
        \beta\left(\divMat{\conductivity}^{(\mathrm{CC})} + \imNum\omega\divMat{\permittivity}^{(\mathrm{CC})}\right) & \beta\divMat{\permittivity}^{(\mathrm{CA})} \\
        \gamma\divMat{\permittivity}^{(\mathrm{AC})} & \gamma\divMat{\permittivity}^{(\mathrm{AA})} \\
    \end{bmatrix}\begin{bmatrix}
        \vectorDOFs^{(\mathrm{C})} \\
        \vectorDOFs^{(\mathrm{A})} \\
    \end{bmatrix}=\mathbf{0}\label{eq:discGenCoul}
\end{equation}
with a weighted divergence matrix
\begin{equation}
    \left(\divMat{\star}\right)_{ij} = \int_{\Omega}\star \DIV\left(\vectorBaseFun_j\right)  \scalarBaseFun_i\,\diffV \label{eq:divMat}
\end{equation}
for $\star\in\{\conductivity,\permittivity\}$.

By utilizing a Lagrange multiplier based approach we can include \eqref{eq:discGenCoul} in the main system \eqref{eq:discrete2} as
\begin{equation}
    \begin{bmatrix}
	    \mathbf{W} &\divMat{}\trans\\
	    \divMat{} & \mathbf{0}
	\end{bmatrix}
	\begin{bmatrix}
	    \vectorDOFs\\
	    \boldsymbol{\lambda} \\
	\end{bmatrix}
	=
	\begin{bmatrix}
	    \mathbf{j}(\scalarDOFs)\\
	    \mathbf{0}
	\end{bmatrix}. \label{eq:langrMulti}
\end{equation}
This increases the size of the linear system and requires the solution of a saddlepoint problem. In order to avoid the use of additional multipliers, we propose a tree-cotree based method to directly include the stabilization in the discrete problem.

\subsection{Tree-Cotree Decomposition}
To address the kernel of $\stiffMatCurl{\reluctivity}$ in all of $\Omega$, we employ a tree-cotree decomposition as described in \cite{Albanese_1988aa,Munteanu_2002aa,Manges_1995aa,Eller_2017aa}. This approach allows us to identify which DOFs remain undetermined and which ones can be determined by only using topological mesh information. The undetermined ones correspond to edges of a spanning tree in the mesh, and we write $\vectorDOFs^{(\mathrm{T})}$ for the respective DOFs. %
The remaining (cotree) DOFs $\vectorDOFs^{(\mathrm{R})}$ must be determined from the magnetostatic problem
\begin{equation}
    \begin{bmatrix}
	    \stiffMatCurl{\reluctivity}^{(\mathrm{RR})} & \stiffMatCurl{\reluctivity}^{(\mathrm{RT})}\\
	    \stiffMatCurl{\reluctivity}^{(\mathrm{TR})} & \stiffMatCurl{\reluctivity}^{(\mathrm{TT})}
	\end{bmatrix}
	\begin{bmatrix}
	    \vectorDOFs^{(\mathrm{R})} \\
	    \vectorDOFs^{(\mathrm{T})}
	\end{bmatrix}
	=
	\begin{bmatrix}
	\mathbf{j}^{(\mathrm{R})}(\mathbf{u})\\
	\mathbf{j}^{(\mathrm{T})}(\mathbf{u})
	\end{bmatrix}. \label{eq:splitMagneto}
\end{equation}
The tree DOFs are often fixed by setting $\vectorDOFs^{(\mathrm{T})}=\mathbf{0}$ \cite{Albanese_1988aa} or by prescribing a relation to the remaining DOFs \cite{Rapetti_2022aa,Munteanu_2002aa}, e.g.,
\begin{equation}
    \vectorDOFs^{(\mathrm{T})}=\stiffMatCurl{\reluctivity}^{(\mathrm{TR})}\left(\stiffMatCurl{\reluctivity}^{(\mathrm{RR})}\right)^{-1}\vectorDOFs^{(\mathrm{R})}.\label{eq:discrOrtho}
\end{equation}
In a second step, \eqref{eq:splitMagneto} is reformulated to compute $\vectorDOFs^{(\mathrm{R})}$. Note that $\stiffMatCurl{\reluctivity}^{(\mathrm{RR})}$ has full rank, and one can show that the second row of \eqref{eq:splitMagneto} is automatically satisfied for a given $\vectorDOFs^{(\mathrm{T})}$ and a $\vectorDOFs^{(\mathrm{R})}$ that is computed from the first line of \eqref{eq:splitMagneto}, see e.g. \cite{Manges_1995aa} or \cite{Rapetti_2022aa}. Note that we construct the tree-cotree splitting for all of $\Omega$ and also employ the following methods in both $\Omega_{\mathrm{C}}$ and $\Omega_{\mathrm{A}}$.

\subsection{Stabilization}
One can reorder the DOFs in \eqref{eq:discrete2} as in \eqref{eq:splitMagneto} to obtain
 \begin{equation}
    \begin{bmatrix}
	    \mathbf{W}^{(\mathrm{RR})} & \mathbf{W}^{(\mathrm{RT})}\\
	    \mathbf{W}^{(\mathrm{TR})} & \mathbf{W}^{(\mathrm{TT})}
	\end{bmatrix}
	\begin{bmatrix}
	    \vectorDOFs^{(\mathrm{R})} \\
	    \vectorDOFs^{(\mathrm{T})}
	\end{bmatrix}
	=
	\begin{bmatrix}
	\mathbf{j}^{(\mathrm{R})}(\mathbf{u})\\
	\mathbf{j}^{(\mathrm{T})}(\mathbf{u})
	\end{bmatrix} \label{eq:TCGsplitSystem}
\end{equation}
and show that $\mathbf{W}^{(\mathrm{RR})}$ has full rank, even for $\omega=0$. The conceptual idea behind the stabilization approach is to replace the  redundant (for $\omega=0$), second line of \eqref{eq:TCGsplitSystem} with the frequency-stable and physically consistent constraint \eqref{eq:discGenCoul}.
Hence, we reorder \eqref{eq:discGenCoul} to arrive at
\begin{equation}
    \begin{bmatrix}
        \divMat{}^{(\mathrm{R})} & \divMat{}^{(\mathrm{T})} \\
    \end{bmatrix}\begin{bmatrix}
	    \vectorDOFs^{(\mathrm{R})}\\
	    \vectorDOFs^{(\mathrm{T})}
	\end{bmatrix}=\mathbf{0}\label{eq:splitCoulomb}
\end{equation}
and combine it with the first line of \eqref{eq:TCGsplitSystem} to obtain
\begin{align}
	\begin{bmatrix}
	    \mathbf{W}^{(\mathrm{RR})} &\mathbf{W}^{(\mathrm{RT})}\\
	    {\divMat{}}^{(\mathrm{R})} & {\divMat{}}^{(\mathrm{T})}
	\end{bmatrix}
	\begin{bmatrix}
	    \mathbf{a}^{(\mathrm{R})}\\
	    \mathbf{a}^{(\mathrm{T})}
	\end{bmatrix}
	=
	\begin{bmatrix}
	    \mathbf{j}^{(\mathrm{R})}(\scalarDOFs) \\
	    \mathbf{0}
	\end{bmatrix}. \label{eq:method1}
\end{align}

\section{Numerical Tests}\label{sec:tests}
To evaluate the performance of the proposed stabilization method, two 3D problems are considered. The first problem is an academic example to investigate the numerical behavior of the method, while the second example deals with a more application-oriented configuration. For our numerical experiments, we employ \texttt{GeoPDEs} \cite{Vazquez_2016aa} and provide the underlying implementations in \cite{Herles_2025ab}. The scaling factors for the scaled divergence matrix \eqref{eq:discGenCoul} are chosen as
{
\begin{equation}
    \beta = 1 + \omega \quad\mathrm{and} \quad \gamma = \left(1 + \omega \right)\frac{\max_{\boldsymbol{x} \in \Omega} \conductivity (\boldsymbol{x}) + \conductivity_\mathrm{art}}{\max_{\boldsymbol{x} \in \Omega } \varepsilon(\boldsymbol{x})},
\end{equation}
with $\conductivity_\mathrm{art} = \SI{1e-6}{}$, to improve the condition of the system.}

\subsection{Academic Example}
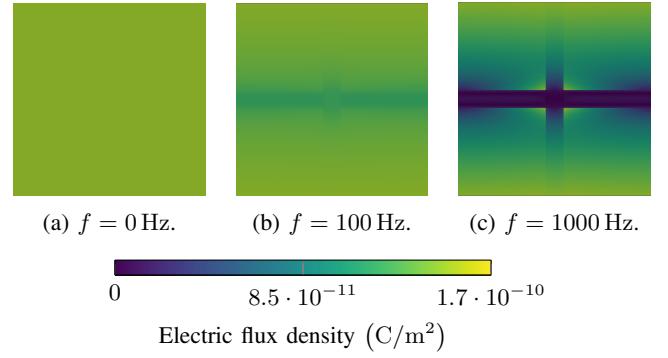
\begin{figure}
    \centering
    \begin{subfigure}[B]{0.32\linewidth}
        \centering
        \includegraphics[width=0.9\linewidth]{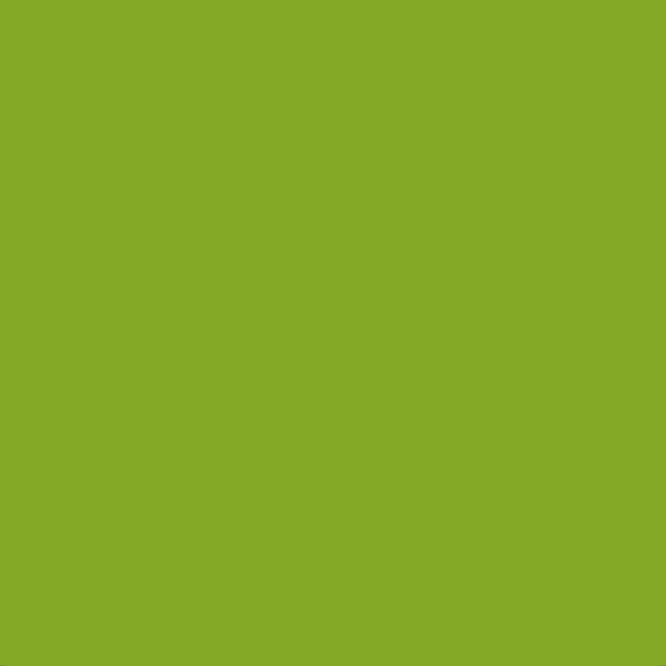}
        \caption{$f=\SI{0}{\hertz}$.}
        \label{fig:Dfield_1}
    \end{subfigure}
    \begin{subfigure}[B]{0.32\linewidth}
        \centering
        \includegraphics[width=0.9\linewidth]{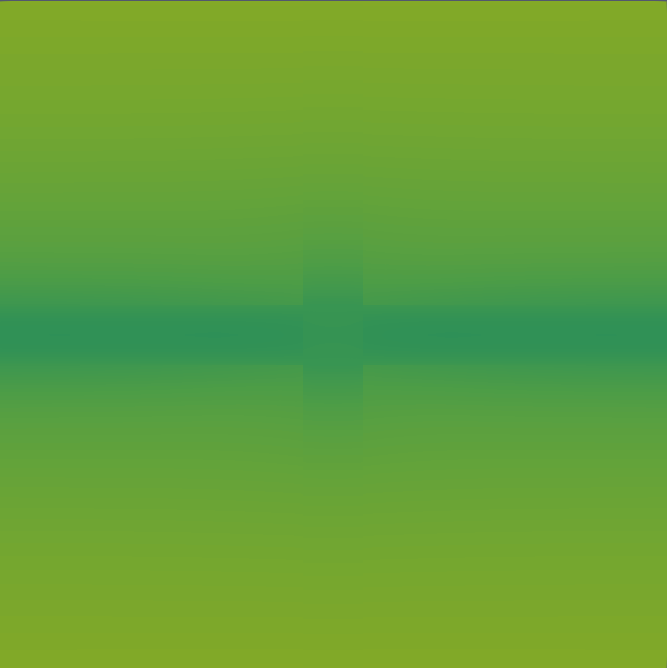}
        \caption{$f=\SI{100}{\hertz}$.}
        \label{fig:Dfield_100}
    \end{subfigure}
    \begin{subfigure}[B]{0.32\linewidth}
        \centering
        \includegraphics[width=0.9\linewidth]{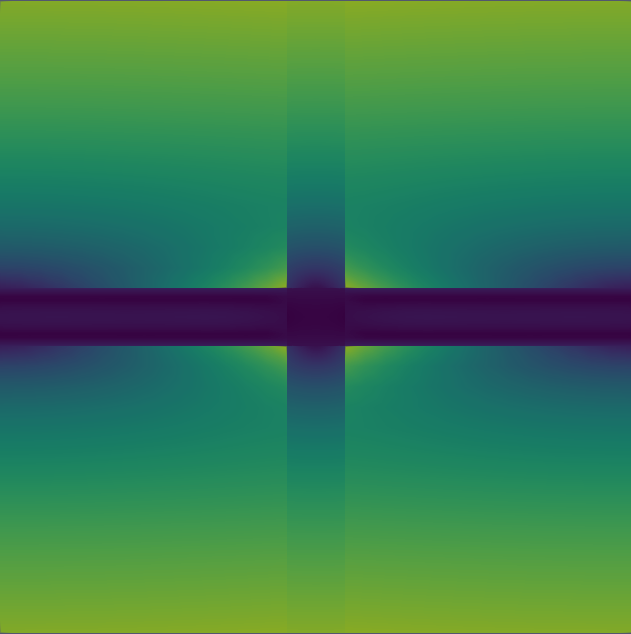}
        \caption{$f=\SI{1000}{\hertz}$.}
        \label{fig:Dfield_1000}
    \end{subfigure}
     \begin{subfigure}[B]{\linewidth}
        \centering
        \begin{tikzpicture}
        \pgfplotsset{every tick label/.append style={font=\small}}
            \begin{axis}[%
  hide axis,
  scale only axis,
  colorbar/width=2mm,
  width = 50mm,
  anchor=south,
  point meta min=0.0,
  point meta max=1.7e-10,
  colormap/viridis,                     %
  colorbar horizontal,                  %
  colorbar sampled,                     %
  colorbar style={
    separate axis lines,
    samples=256,                        %
    xlabel=\small{Electric flux density $\big(\SI{}{\coulomb/\meter^{2}}\big)$},
    xtick={0.0, 0.85e-10, 1.7e-10},
    scaled x ticks=false
  },
]
  \addplot [draw=none] coordinates {(0,0)};
\end{axis}
        \end{tikzpicture}
    \end{subfigure}
  \caption{Electric flux density $\norm{\Dfield}{2}$ for different excitation frequencies.}
\end{figure}

The first test problem consists of three conducting bars in a dielectric box. The configuration is shown in \autoref{fig:academicTestEx} and is based on the test problem of \cite{Ostrowski_2021aa}. For the magnetic vector potential, we employ a homogeneous Dirichlet boundary $\Afield \times \nfield = \zerofield$ on the whole boundary $\partial\Omega$. It is constructed such that the flux density $\Dfield_{\mathrm{e}}$  does not depend on $\omega$ and is constant in space. For a more detailed derivation, see \cite{Ostrowski_2021aa}. {The problem is discretized using 27 patches and {6084 third order basis functions}.} \autoref{fig:Dfield_1} shows $\norm{\Dfield}{2}$ for $f=\SI{0}{\hertz}$ where no frequency-dependent contributions from $\Dfield_{\mathrm{m}}$ are present such that only the constant $\Dfield_{\mathrm{e}}$ contributes. In \autoref{fig:Dfield_100} for $f=\SI{100}{\hertz}$ and \autoref{fig:Dfield_1000} for $f=\SI{1000}{\hertz}$, one can see how the resulting flux density is being pulled out of the conductor due to the correction~$\Dfield_{\mathrm{m}}$.

\def\LEN{5}
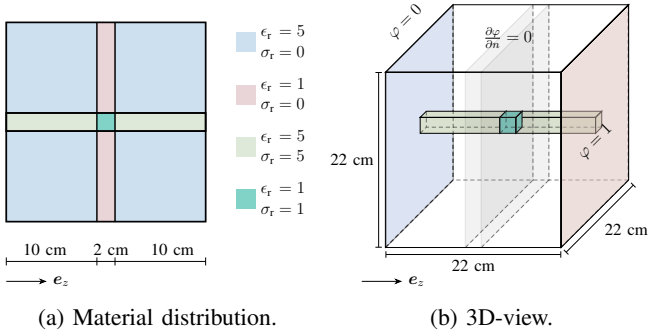
\begin{figure}
    \centering
      \subfloat[Material distribution.\label{fig:Cube_mat}]{%
       \scalebox{0.53}{
                    \begin{tikzpicture}
                        \tikzset{every node/.append style = {font=\large}}
                		\draw[black, very thick] (0,0) rectangle (\LEN,\LEN);

                		\draw (0,-1) -- (\LEN,-1);
                		\node[xshift={3mm}] at (\LEN*2/15,-0.7) {10 cm};
                		\node[xshift={-1mm}] at (\LEN*6/11,-0.7) {2 cm};
                		\node[xshift={-5mm}] at (\LEN*14/15,-0.7) {10 cm};
                		\draw (0,-1.1) -- (0,-0.9);
                		\draw (\LEN*10/22,-1.1) -- (\LEN*10/22,-0.9);
                		\draw (\LEN*12/22,-1.1) -- (\LEN*12/22,-0.9);
                		\draw (\LEN,-1.1) -- (\LEN,-0.9);

                		\draw[fill=TUDa-1b,fill opacity=0.2] (0,0) rectangle (\LEN*10/22,\LEN*10/22);
                        \draw[fill=TUDa-1b,fill opacity=0.2] (\LEN*12/22,0) rectangle (\LEN,\LEN*10/22);
                        \draw[fill=TUDa-1b,fill opacity=0.2] (0,\LEN*12/22) rectangle (\LEN*10/22,\LEN);
                        \draw[fill=TUDa-1b,fill opacity=0.2] (\LEN*12/22,\LEN*12/22) rectangle (\LEN,\LEN);

                        \draw[fill=TUDa-9d, fill opacity=0.2] (\LEN*10/22,0) rectangle (\LEN*12/22,\LEN*10/22);
                        \draw[fill=TUDa-9d, fill opacity=0.2] (\LEN*10/22,\LEN*12/22) rectangle (\LEN*12/22,\LEN);

                		\draw[fill=TUDa-4d, fill opacity=0.2] (0,\LEN*10/22) rectangle (\LEN*10/22,\LEN*12/22);
                		\draw[fill=TUDa-3b, fill opacity=0.5] (\LEN*10/22,\LEN*10/22) rectangle (\LEN*12/22,\LEN*12/22);
                		\draw[fill=TUDa-4d, fill opacity=0.2] (\LEN*12/22,\LEN*10/22) rectangle (\LEN,\LEN*12/22);

                        \node[fill=TUDa-1b, rectangle, inner sep= 0pt, minimum width=5mm, minimum height=5mm, fill opacity=0.2](l1) at (\LEN+1,\LEN-0.5){};
                        \node[fill=TUDa-9d, rectangle, inner sep= 0pt, minimum width=5mm, minimum height=5mm, fill opacity=0.2,anchor=north](l2) at ($(l1.south)-(0,0.8)$){};
                        \node[fill=TUDa-4d, rectangle, inner sep= 0pt, minimum width=5mm, minimum height=5mm, fill opacity=0.2,anchor=north](l3) at ($(l2.south)-(0,0.8)$){};
                        \node[fill=TUDa-3b, rectangle, inner sep= 0pt, minimum width=5mm, minimum height=5mm, fill opacity=0.5,anchor=north](l4) at ($(l3.south)-(0,0.8)$){};

                        \node[anchor=west] at (l1.north east){$\permittivity_\text{r}\,=5$};
                        \node[anchor=west] at (l1.south east){$\conductivity_\text{r}=0$};
                        \node[anchor=west] at (l2.north east){$\permittivity_\text{r}\,=1$};
                        \node[anchor=west] at (l2.south east){$\conductivity_\text{r}=0$};
                        \node[anchor=west] at (l3.north east){$\permittivity_\text{r}\,=5$};
                        \node[anchor=west] at (l3.south east){$\conductivity_\text{r}=5$};
                        \node[anchor=west] at (l4.north east){$\permittivity_\text{r}\,=1$};
                        \node[anchor=west] at (l4.south east){$\conductivity_\text{r}=1$};

                        \draw[black, very thick] (0,\LEN*10/22) rectangle (\LEN,\LEN*12/22);

                		\draw[-{Stealth[black]}] (0,-1.5)   -- node[pos=1.35, xshift={-3mm}, anchor=west]{$\boldsymbol{e}_z$} (1,-1.5);
                	\end{tikzpicture}
                }
       }
    \hfill
  \subfloat[3D-view.\label{fig:Cube3D}]{%
       \scalebox{0.53}{%
                    \begin{tikzpicture}
                        \tikzset{every node/.append style = {font=\large}}

                        \pic [fill=TUDa-1b!25, text=TUDa-1b, draw=black,opacity=0.5] at (-4.4,0) {cuboid={width=0.1, height=44, depth=44}};
                		\pic [fill=TUDa-4d!20, text=TUDa-1b, draw=black] at (+0.866,-2+0.866) {cuboid={width=20, height=4, depth=4}};
                		\pic [fill=TUDa-4d!20, text=TUDa-1b, draw=black] at (+0.866-2-0.4,-2+0.866) {cuboid={width=20, height=4, depth=4}};
                		\pic [fill=none, text=black, draw=black] at (0,0) {annotated cuboid={width=44, height=44, depth=44}};
                		\pic [fill=TUDa-3b!50, text=TUDa-1b, draw=black] at (+0.866-2,-2+0.866) {cuboid={width=4, height=4, depth=4}};

                		\pic [fill=TUDa-9d!25, text=TUDa-1b, draw=black,opacity=0.5] at (0,0) {cuboid={width=0.1, height=44, depth=44}};

                		\pic [fill=gray, text=TUDa-1b, draw=black, opacity=0.1] at (-2.4,0) {cuboid={width=0, height=44, depth=44}};
                		\pic [fill=gray, text=TUDa-1b, draw=black, opacity=0.1] at (-2,0) {cuboid={width=0, height=44, depth=44}};

                		\draw[-{Stealth[black]}] (0-\LEN,-0.25-\LEN)   -- node[pos=1.35, xshift={-3mm}, anchor=west]{$\boldsymbol{e}_z$} (1-\LEN,-0.25-\LEN);

                		\node[xshift={14mm},yshift={29mm},rotate=45] at (-\LEN-0.3,-\LEN*0.5+0.866) {$\phifield = 0$};
                		\node[xshift={-2mm},rotate=45] at (+0.866 + 0.2,-\LEN*0.5+0.866-0.1) {$\phifield = 1$};
                		\node[xshift={3mm}] at (-\LEN*0.5+0.866,+0.866) {$\frac{\partial \phifield}{\partial \nfield} = 0$};
                	\end{tikzpicture}
                }
       }
  \caption{Academic test example problem setup.}
  \label{fig:academicTestEx}
\end{figure}

\begin{figure}
    \centering
 \pgfplotstableread[col sep=comma]{data/div_cond_residuals_cube_p2_s2.csv}\residual
                \begin{tikzpicture}
                	\begin{loglogaxis}[tudalineplot, width=1\linewidth, height=12em, xlabel={\scriptsize$f$ in Hz}, ylabel={\scriptsize Residual $\delta_{\divMat{}}$}, xticklabel style={font=\scriptsize}, xmin=5.5e-7,xmax=2e14, ymin=1e-23, ymax=1e4, xtick distance=10^(4), ytick distance=1e4, xlabel shift={-1mm}, ylabel shift={-1mm},legend style={at={(1,0.85)},anchor=east}, legend style={nodes={scale=0.7, transform shape}}]

                 \addplot+[mark size=1.5mm,mark=*,mark options={fill=TUDa-8a,scale=0.3}, TUDa-8a, line width=1pt] table[x index = 0, y index = 7] {\residual};

            \addplot+[mark size=1.5mm,mark=square*,mark options={fill=TUDa-4a,scale=0.3}, TUDa-4a, line width=1pt] table[x index = 0, y index = 9] {\residual};

            \addlegendentry{original};
            \addlegendentry{stabilized};
        \end{loglogaxis}
    \end{tikzpicture}
    \caption{Residual of discrete implicit constraint of the academic test example over frequency.}
    \label{fig:div}
\end{figure}
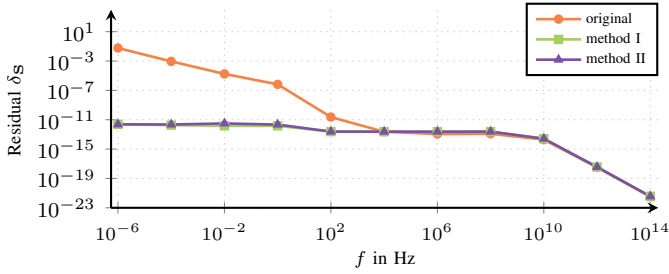

{In the following, ``original'' or ``original system'' refer to \eqref{eq:discrete2}.} \autoref{fig:condNums} shows the condition number of the original system matrix {and of the stabilized method.
For low frequencies, the original system matrix {is} ill-conditioned {or even singular}, while the condition number of the stabilized method {is relatively small} for $f \rightarrow \SI{0}{Hz}$. In the limit $f = \SI{0}{Hz}$, we obtain a condition number of approximately $\SI{4.87e5}{}$, while the original system matrix is singular. For high frequencies, the condition number of the stabilized system increases slightly in comparison to the one of the original.} This effect is well known, see for example \cite{Hiptmair_2000aa} and \cite{Munteanu_2002aa}. {As another verification,} \autoref{fig:div} shows the residual of the discretely evaluated generalized Coulomb gauge \eqref{eq:strongGenCoul} namely
\begin{equation*}
    \delta_{\divMat{}} = \left\Vert \divMat{\complexConductivity}\mathbf{a}\right\Vert_{2}\overset{!}{=}0,
\end{equation*}
where $\left\Vert\cdot\right\Vert_2$ denotes the $\ell_2$-norm.
The implicit constraint is well satisfied over the whole frequency range for { the stabilized method ($\delta_{\divMat{}}<\SI{1e-11}{}$). Even in the static limit $f=\SI{0}{Hz}$, {the stabilized system} remains consistent ($\delta_{\divMat{}}=\SI{1.52e-12}{}$).} On the other hand, the residual of the original formulation grows significantly for small frequencies up to $\delta_{\divMat{}}\geq\SI{1e-2}{}$.

\begin{figure}
    \centering
 \pgfplotstableread[col sep=comma]{data/cond_fullyScaled.csv}\cond
                \begin{tikzpicture}
                	\begin{loglogaxis}[tudalineplot, width=1\linewidth, height=12em, xlabel={\scriptsize$f$ in Hz}, ylabel={\scriptsize Cond. Number}, xticklabel style={font=\scriptsize}, xmin=5.5e-7,xmax=2e14, ymin=5e3, ymax=5e20, xtick distance=10^(4), ytick distance=10^(4), xlabel shift={-1mm}, ylabel shift={-1mm},legend style={at={(1,0.85)},anchor=east}, legend style={nodes={scale=0.7, transform shape}}]

                 \addplot+[mark size=1.5mm,mark=*,mark options={fill=TUDa-8a,scale=0.3}, TUDa-8a, line width=1pt] table[x index = 0, y index = 1] {\cond};

                        \addplot+[mark size=1.5mm,mark=square*,mark options={fill=TUDa-4a,scale=0.3}, TUDa-4a, line width=1pt] table[x index = 0, y index = 7] {\cond};

                    \addlegendentry{original};
                		\addlegendentry{method I};
                	\end{loglogaxis}
                \end{tikzpicture}

                    \caption{Condition number of different system matrices of academic test example over frequency.}
    \label{fig:condNums}
\end{figure}
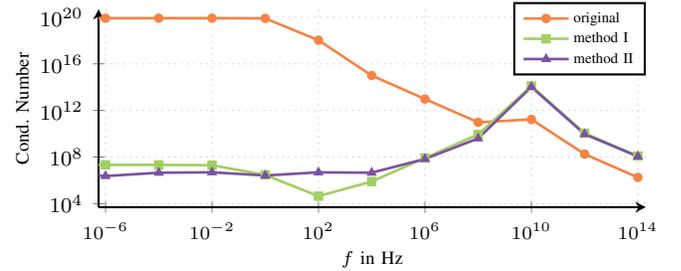

{
\subsection{Convergence Study for Constructed Example}
In this subsection, we verify whether or not the stabilized method computes correct numerical solutions to the PDE at hand. Let $\Omega=\left(\frac{\pi}{2},\frac{3\pi}{2}\right)^3$ be given as our domain and we assume
\begin{align*}
    \Afield_{\mathrm{ana}}&=\begin{bmatrix}
        \sin(x)\cos(y)\cos(z) \\
        -2\cos(x)\sin(y)\cos(z) \\
        \cos(x)\cos(y)\sin(z) \\
    \end{bmatrix} \\
    \phifield_{\mathrm{ana}}&=\cos(x)\cos(y)\cos(z)
\end{align*}
for the solution prescribed on $\Omega$. Both $\Afield_{\mathrm{ana}}\times\nfield=\zerofield$ and $\phifield_{\mathrm{ana}}=0$, are then statisfied on $\partial\Omega$. Furthermore, $\DIV\left(\complexConductivity\Afield_{\mathrm{ana}}\right)=0$ as long as $\complexConductivity=\mathrm{const.}$ is chosen in the setup. We choose $\reluctivity=\frac{1}{\mu_0}$, $\permittivity=\permittivity_0$ and use {either $\conductivity=0$ or $\conductivity=\SI{6e7}{}$.} This allows for the closed-form computation of the source terms $\Jfield^{\mathrm{s}}$ and $\rhofield^{\mathrm{s}}$ such that a numerical simulation converges to the corresponding analytical solution.
}

\begin{figure}
    \centering
    \def\errIDOne{4}
    \def\errIDTwo{5}
    \def\errIDThree{6}
    \begin{subfigure}[T]{0.49\linewidth}
        \centering
        \begin{tikzpicture}
			\begin{loglogaxis}[
			    xmin=1.5, xmax=20, ymin=5e-5, ymax=100, log origin=infty, xlabel={\small$s_{\mathrm{h}}$}, ylabel={\small Error}, xtick={2,4,8,16}, xticklabels={2,4,8,16}, width=0.99\linewidth, height=4cm
            ]

            \addplot+[solid, mark=square, blue] table[x index = 0,y index = \errIDOne, col sep=comma] {data/errors_original_sig0_f1000000_mat.csv};\label{plot:p1}
			\addplot+[solid, mark=square, red] table[x index = 0,y index = \errIDTwo, col sep=comma] {data/errors_original_sig0_f1000000_mat.csv};\label{plot:p2}
			\addplot+[solid, mark=square, brown] table[x index = 0,y index = \errIDThree, col sep=comma] {data/errors_original_sig0_f1000000_mat.csv};\label{plot:p3}

            \addplot+[solid, mark=x, blue] table[x index = 0,y index = \errIDOne, col sep=comma] {data/errors_stabilization_sig0_f1000000_mat.csv}; \label{plot:p4}
			\addplot+[solid, mark=x, red] table[x index = 0,y index = \errIDTwo, col sep=comma] {data/errors_stabilization_sig0_f1000000_mat.csv}; \label{plot:p5}
			\addplot+[solid, mark=x, brown] table[x index = 0,y index = \errIDThree, col sep=comma] {data/errors_stabilization_sig0_f1000000_mat.csv}; \label{plot:p6}

            \addplot+[black, dashed, mark=none, samples at={2,3,4,6,8,16}] {10^(0.8-log10(x))};\label{plot:p7}
            \addplot+[black, dashed, mark=none, samples at={2,3,4,6,8,16}, forget plot] {10^(0.8-2*log10(x))};
            \addplot+[black, dashed, mark=none, samples at={2,3,4,6,8,16}, forget plot] {10^(0.8-3*log10(x))};

			\end{loglogaxis}
        \end{tikzpicture}
        \centering
        \caption{$\conductivity = 0$, $f=\SI{1}{\mega\hertz}$}
        \label{fig:conv1}
    \end{subfigure}
    \begin{subfigure}[T]{0.49\linewidth}
        \centering
        \begin{tikzpicture}
			\begin{loglogaxis}[
			    xmin=1.5, xmax=20, ymin=5e-5, ymax=1e2, log origin=infty, xlabel={\small$s_{\mathrm{h}}$}, ylabel={\small Error}, xtick={2,4,8,16}, xticklabels={2,4,8,16}, width=0.99\linewidth, height=4cm,
            ]

            \addplot+[solid, mark=square, blue] table[x index = 0,y index = \errIDOne, col sep=comma] {data/errors_original_sig0_f10_mat.csv};
			\addplot+[solid, mark=square, red] table[x index = 0,y index = \errIDTwo, col sep=comma] {data/errors_original_sig0_f10_mat.csv};
			\addplot+[solid, mark=square, brown] table[x index = 0,y index = \errIDThree, col sep=comma] {data/errors_original_sig0_f10_mat.csv};

            \addplot+[solid, mark=x, blue] table[x index = 0,y index = \errIDOne, col sep=comma] {data/errors_stabilization_sig0_f10_mat.csv};
			\addplot+[solid, mark=x, red] table[x index = 0,y index = \errIDTwo, col sep=comma] {data/errors_stabilization_sig0_f10_mat.csv};
			\addplot+[solid, mark=x, brown] table[x index = 0,y index = \errIDThree, col sep=comma] {data/errors_stabilization_sig0_f10_mat.csv};

            \addplot+[black, dashed, mark=none, samples at={2,3,4,6,8,16}] {10^(0.8-log10(x))};
            \addplot+[black, dashed, mark=none, samples at={2,3,4,6,8,16}, forget plot] {10^(0.8-2*log10(x))};
            \addplot+[black, dashed, mark=none, samples at={2,3,4,6,8,16}, forget plot] {10^(0.8-3*log10(x))};

			\end{loglogaxis}
        \end{tikzpicture}
        \centering
        \caption{$\conductivity = 0$, $f=\SI{10}{\hertz}$}
        \label{fig:conv2}
    \end{subfigure}
    \begin{subfigure}[T]{0.49\linewidth}
        \centering
        \begin{tikzpicture}
			\begin{loglogaxis}[
			    xmin=1.5, xmax=20, ymin=5e-5, ymax=100, log origin=infty, xlabel={\small$s_{\mathrm{h}}$}, ylabel={\small Error}, xtick={2,4,8,16}, xticklabels={2,4,8,16}, width=0.99\linewidth, height=4cm
            ]

            \addplot+[solid, mark=square, blue] table[x index = 0,y index = \errIDOne, col sep=comma] {data/errors_original_sig6e7_f10_mat.csv};
			\addplot+[solid, mark=square, red] table[x index = 0,y index = \errIDTwo, col sep=comma] {data/errors_original_sig6e7_f10_mat.csv};
			\addplot+[solid, mark=square, brown] table[x index = 0,y index = \errIDThree, col sep=comma] {data/errors_original_sig6e7_f10_mat.csv};

            \addplot+[solid, mark=x, blue] table[x index = 0,y index = \errIDOne, col sep=comma] {data/errors_stabilization_sig6e7_f10_mat.csv};
			\addplot+[solid, mark=x, red] table[x index = 0,y index = \errIDTwo, col sep=comma] {data/errors_stabilization_sig6e7_f10_mat.csv};
			\addplot+[solid, mark=x, brown] table[x index = 0,y index = \errIDThree, col sep=comma] {data/errors_stabilization_sig6e7_f10_mat.csv};

            \addplot+[black, dashed, mark=none, samples at={2,3,4,6,8,16}] {10^(0.8-log10(x))};
            \addplot+[black, dashed, mark=none, samples at={2,3,4,6,8,16}, forget plot] {10^(0.8-2*log10(x))};
            \addplot+[black, dashed, mark=none, samples at={2,3,4,6,8,16}, forget plot] {10^(0.8-3*log10(x))};

			\end{loglogaxis}
        \end{tikzpicture}
        \centering
        \caption{$\conductivity = \SI{6e7}{}$, $f=\SI{10}{\hertz}$}
        \label{fig:conv3}
    \end{subfigure}
    \begin{subfigure}[T]{0.49\linewidth}
    \vspace{2mm}
    \centering
    \begin{tikzpicture}
        \matrix[
		matrix of nodes,
		draw, font=\small,
		inner sep=0.2em,anchor=south west,
		ampersand replacement={\&}, column 1/.style={anchor=east}, column 2/.style = {anchor=west} ]at(0,-2)
					{orig: \hspace{0.3cm}
                    \& \ref*{plot:p1} \hspace{0.25cm} $p=1$ \\
					\& \ref*{plot:p2} \hspace{0.25cm} $p=2$ \\
                    \& \ref*{plot:p3} \hspace{0.25cm} $p=3$ \\[1mm]
				    stab.: \hspace{0.3cm}
                    \& \ref*{plot:p4} \hspace{0.25cm} $p=1$ \\
					\& \ref*{plot:p5} \hspace{0.25cm} $p=2$ \\
                    \& \ref*{plot:p6} \hspace{0.25cm} $p=3$ \\[1mm]
                    \& \ref*{plot:p7} \hspace{0.25cm} $\mathcal{O}\left(s_{\mathrm{h}}^{-p}\right)$ \\};    \end{tikzpicture}
    \end{subfigure}
    \caption{$H(\Omega,\mathrm{curl})$-Error for different conductivities, frequencies and degrees over mesh refinement $s_{\mathrm{h}}\propto h^{-1}$.}
    \label{fig:convergence}
\end{figure}

{\autoref{fig:convergence} shows the $H(\Omega,\mathrm{curl})$-Error of the numerically computed magnetic vector potential using the original and stabilized formulations over the number of subdivisions $s_h$ for different degrees $p$ and problem setups. One can see in \autoref{fig:conv1} that both formulations remain stable even in non-conducting domains and converge optimally for a frequency of $f=\SI{1}{\mega\hertz}$. For lower frequencies, like $f=\SI{10}{Hz}$, the original formulation starts to break down while the stabilized formulation preserves its optimal convergence behavior {as shown in \autoref{fig:conv2}. Note that the error for $p=3$ completely diverges and is not contained in the visualized error region anymore}. \autoref{fig:conv3} shows that the breakdown occurs for much lower frequencies in conducting domains and that both formulations are stable for $f=\SI{10}{Hz}$.
}

\subsection{Simulation of Setup with Planar Coil} \label{sec:numPlanarCoil}
\begin{figure}
    \centering
    \begin{tikzpicture}
        \node[inner sep=0pt] (test) at (0,0) {\includegraphics[width=\linewidth]{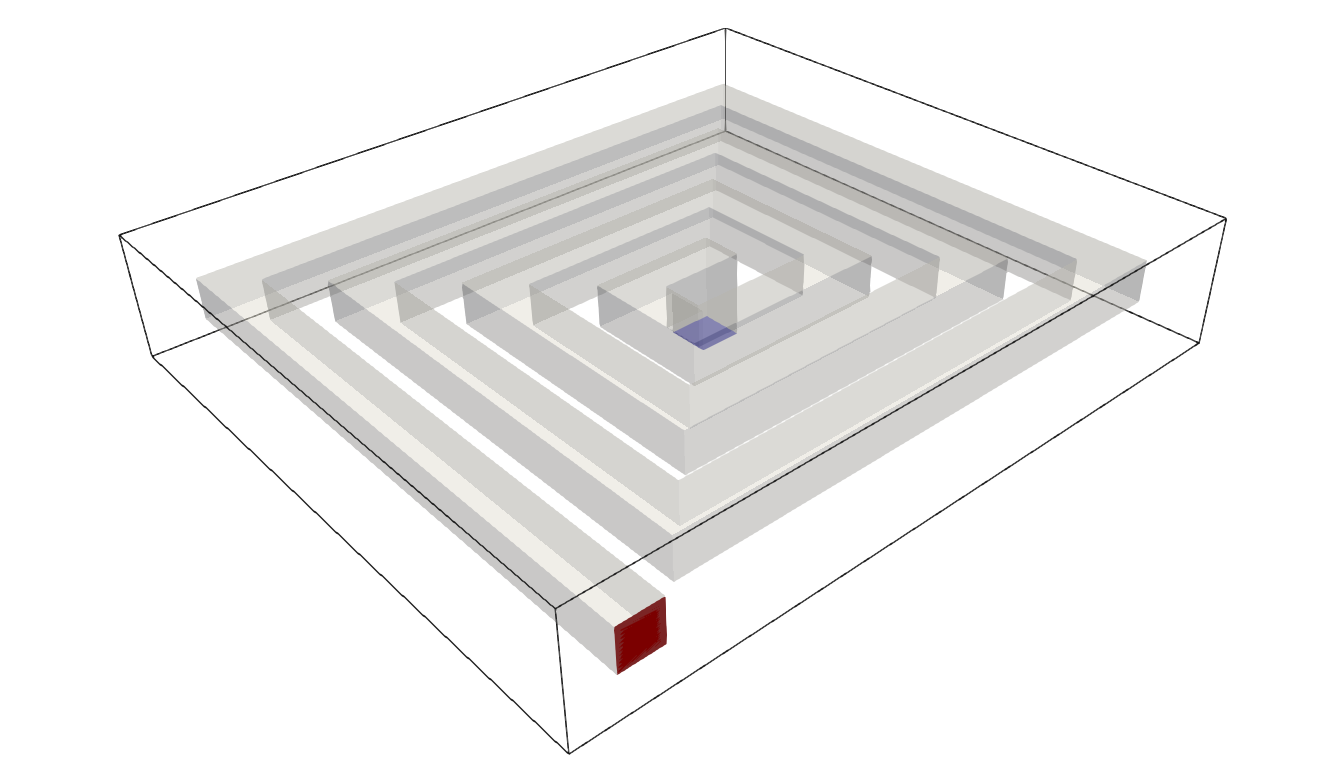}};
        \node[inner sep=0pt] (test) at (-0.2,0) {$\Gamma_\mathrm{G}$};
        \node[inner sep=0pt] (test) at (0.4,-1.5) {$\Gamma_\mathrm{E}$};
    \end{tikzpicture}
    \caption{Problem setup of the planar coil.}
    \label{fig:planarCoil}
\end{figure}

As a more application-motivated example, we consider a test case based on a planar coil. The problem setup is shown in \autoref{fig:planarCoil}. The copper coil has a square cross-section of $\SI{3}{mm} \times \SI{3}{mm}$ and a conductivity of $\conductivity = \SI{6e7}{\frac{S}{m}}$. It has three
{turns}
and is surrounded by air. {The airbox has dimensions $\SI{45}{mm} \times \SI{51}{mm} \times \SI{9}{mm}$.} A voltage is applied to the coil, resulting in boundary conditions $\phifield=\SI{0}{V}$ on $\Gamma_{\mathrm{G}}$ and $\phifield=(\sqrt{2})^{-1}\phi_{\mathrm{amp}}$ on $\Gamma_\mathrm{E}$ with $\phi_{\mathrm{amp}}=\SI{1}{V}$ oscillating with frequency  $f$. On the remaining boundary $\partial\Omega\setminus \left(\Gamma_\mathrm{G}\cup\Gamma_\mathrm{E}\right)$, we set $\GRAD\phifield\cdot\nfield=0$. For the magnetic vector potential, we assume perfect magnetic boundary conditions, i.e., $\CURL\Afield\times\nfield=\zerofield$. We do not impress a source current density. {The problem is discretized using 144 patches and {590,326 second order basis functions.}}

\begin{figure}
    \vspace*{-2em}
    \centering
    \begin{subfigure}[B]{0.49\linewidth}
        \centering
        \includegraphics[width=0.9\linewidth,trim={2.8cm 3.5cm 3.5cm 3.5cm},clip]{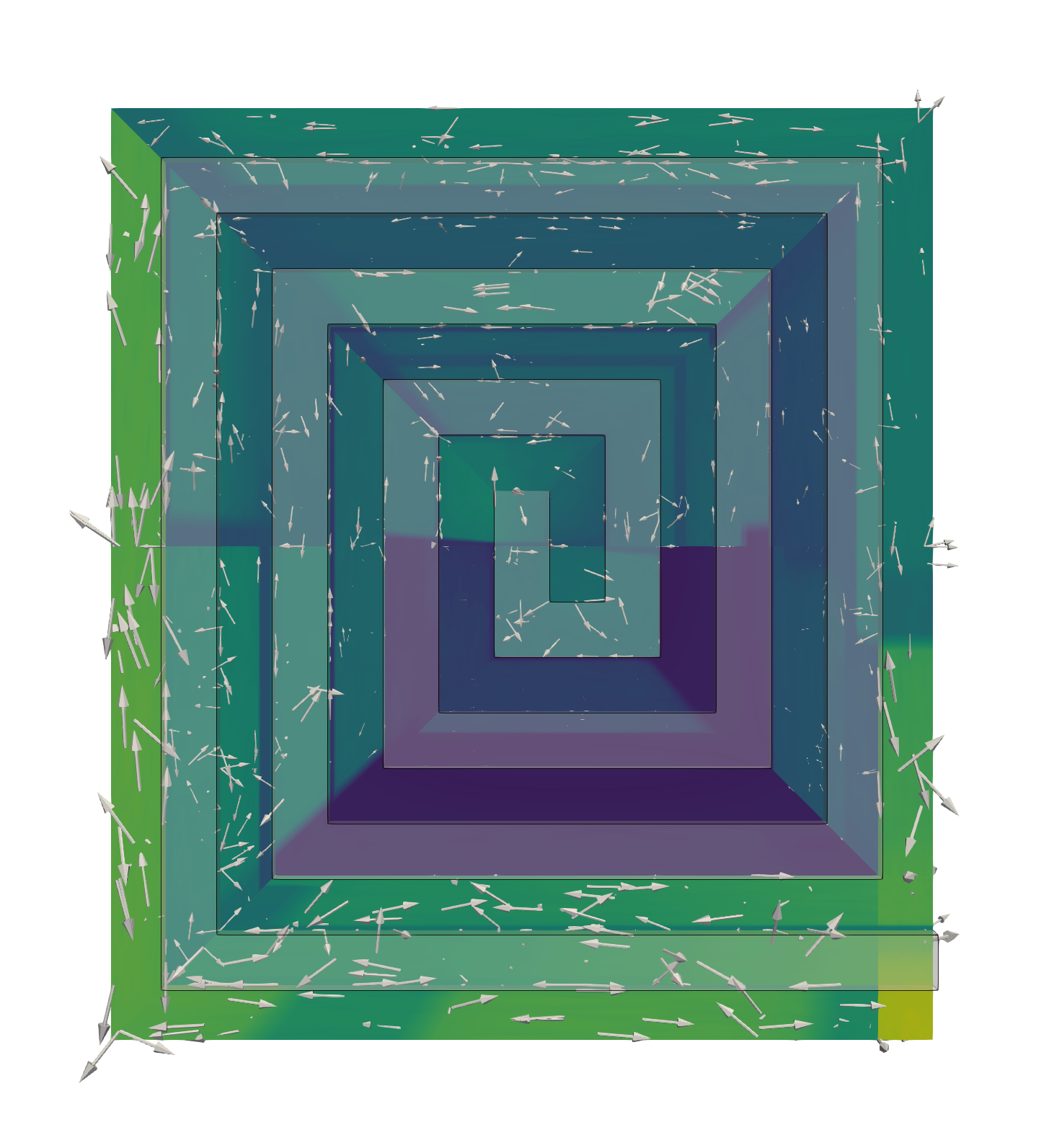}
        \\
    \begin{tikzpicture}[every node/.style={scale=0.8}]
        \pgfplotsset{every tick label/.append style={font=\small}}
            \begin{axis}[%
  hide axis,
  scale only axis,
  colorbar/width=2mm,
  width = 33mm,
  anchor=south,
  point meta min=0.0,
  point meta max=1,
  colormap/viridis,                     %
  colorbar horizontal,                  %
  colorbar sampled,                     %
  colorbar style={
    separate axis lines,
    samples=256,                        %
    xlabel=\small{Magnetic flux density $\big(\SI{}{\tesla}\big)$},label style={font=\tiny},
    xtick={0.0,0.5, 1},
    xticklabels={%
        $\num{1e-10}$,
        $\num{1e30}$,
        $\num{1e70}$},
    scaled x ticks=false
  },
]
  \addplot [draw=none] coordinates {(0,0)};
\end{axis}
        \end{tikzpicture}
        \caption{$\norm{\Bfield}{2}$, non-stabilized.}
        \label{fig:Bfield_pC_unstab}
    \end{subfigure}
     \begin{subfigure}[B]{0.49\linewidth}
        \centering
        \includegraphics[width=0.9\linewidth,trim={2.8cm 3.5cm 3.5cm 3.5cm},clip]{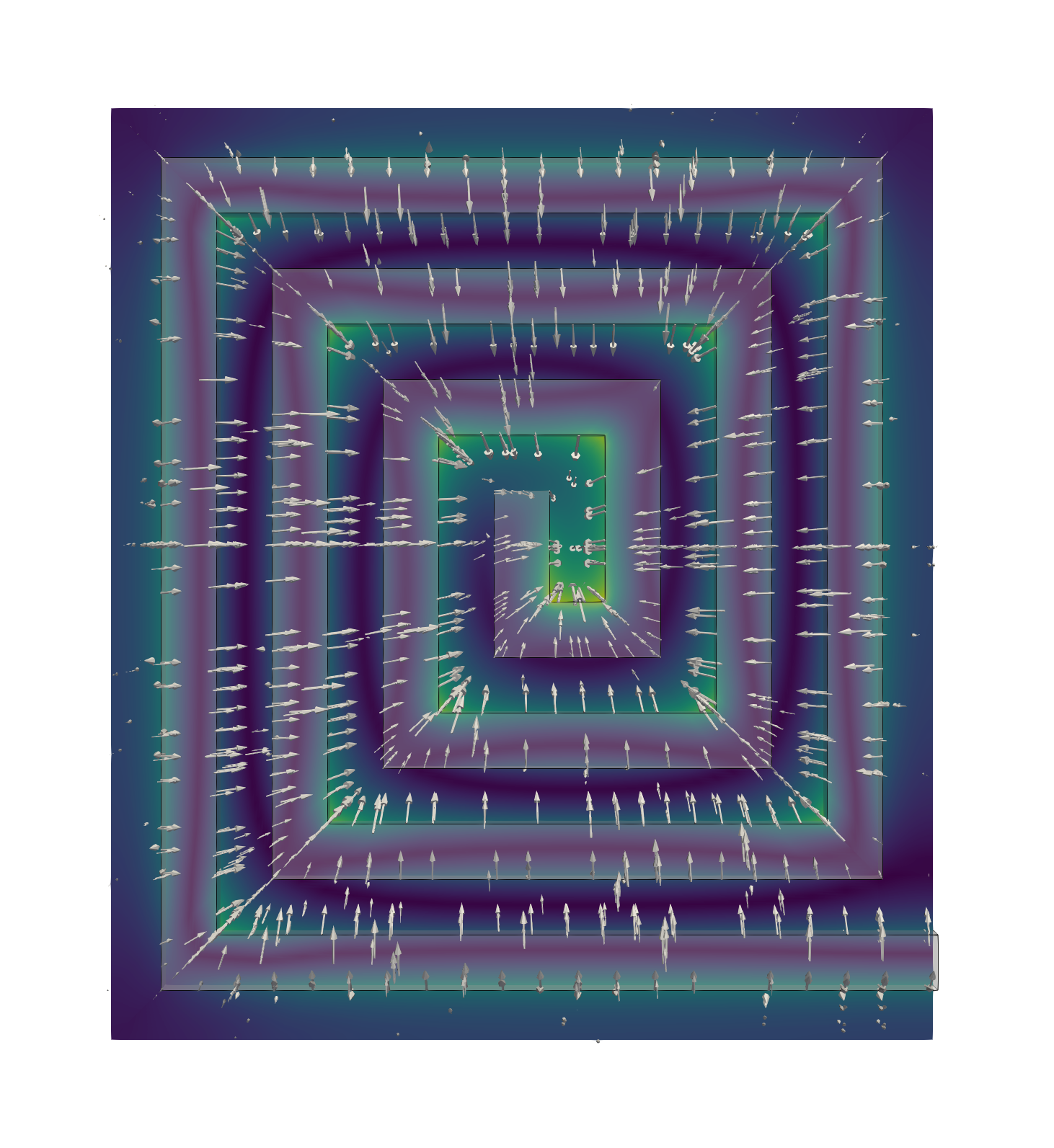}\\
        \hspace{1mm}
    \begin{tikzpicture}[every node/.style={scale=0.8}]
        \pgfplotsset{every tick label/.append style={font=\small}}
            \begin{axis}[%
  hide axis,
  scale only axis,
  colorbar/width=2mm,
  width = 33mm,
  anchor=south,
  point meta min=0.0,
  point meta max=0.42,
  colormap/viridis,                     %
  colorbar horizontal,                  %
  colorbar sampled,                     %
  colorbar style={
    separate axis lines,
    samples=256,                        %
    xlabel=\small{Magnetic flux density $\big(\SI{}{\tesla}\big)$},
    xtick={0.0,0.21, 0.42},
    scaled x ticks=false
  },
]
  \addplot [draw=none] coordinates {(0,0)};
\end{axis}
        \end{tikzpicture}
        \caption{$\norm{\Bfield}{2}$, stabilized.}
        \label{fig:Bfield_pC}
    \end{subfigure} \\[1em]
    \begin{subfigure}[B]{0.49\linewidth}
        \centering
        \includegraphics[width=0.9\linewidth,trim={2.8cm 3.5cm 3.5cm 3.5cm},clip]{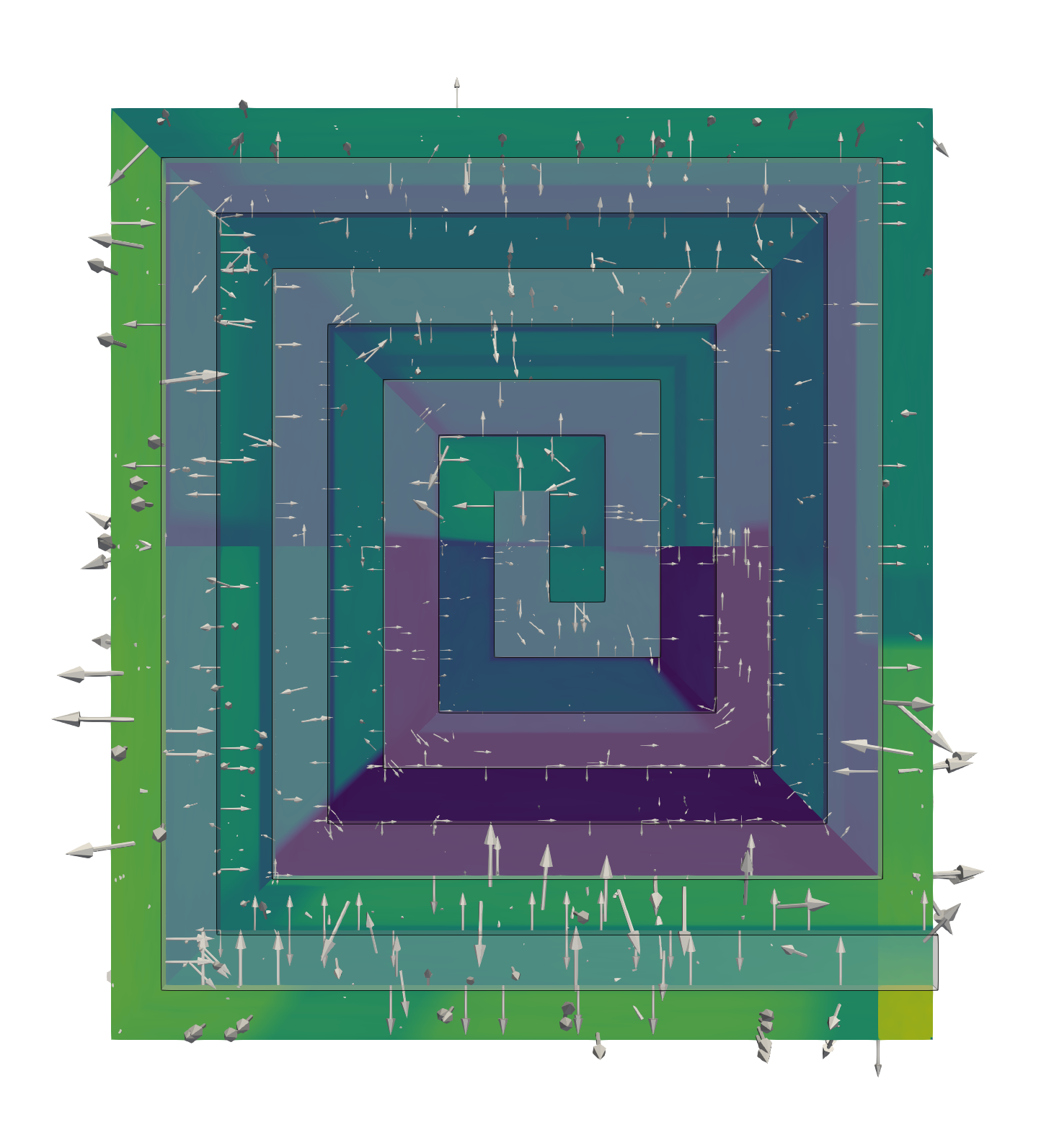}
        \begin{tikzpicture}[every node/.style={scale=0.8}]
        \pgfplotsset{every tick label/.append style={font=\small}}
            \begin{axis}[%
  hide axis,
  scale only axis,
  colorbar/width=2mm,
  width = 33mm,
  anchor=south,
  point meta min=0.0,
  point meta max=1,
  colormap/viridis,                     %
  colorbar horizontal,                  %
  colorbar sampled,                     %
  colorbar style={
    separate axis lines,
    samples=256,                        %
    xlabel={\small Electric field strength $\big(\SI{}{\volt/\meter}\big)$},
    xtick={0.0,0.5, 1},
    xticklabels={%
        $\num{1e-15}$,
        $\num{1e30}$,
        $\num{1e75}$},
    scaled x ticks=false
  },
]
  \addplot [draw=none] coordinates {(0,0)};
\end{axis}
        \end{tikzpicture}
        \caption{$\norm{\Efield}{2}$, non-stabilized.}
        \label{fig:Efield_pC_unstab}
    \end{subfigure}
    \begin{subfigure}[B]{0.49\linewidth}
        \centering
        \includegraphics[width=0.9\linewidth,trim={2.8cm 3.5cm 3.5cm 3.5cm},clip]{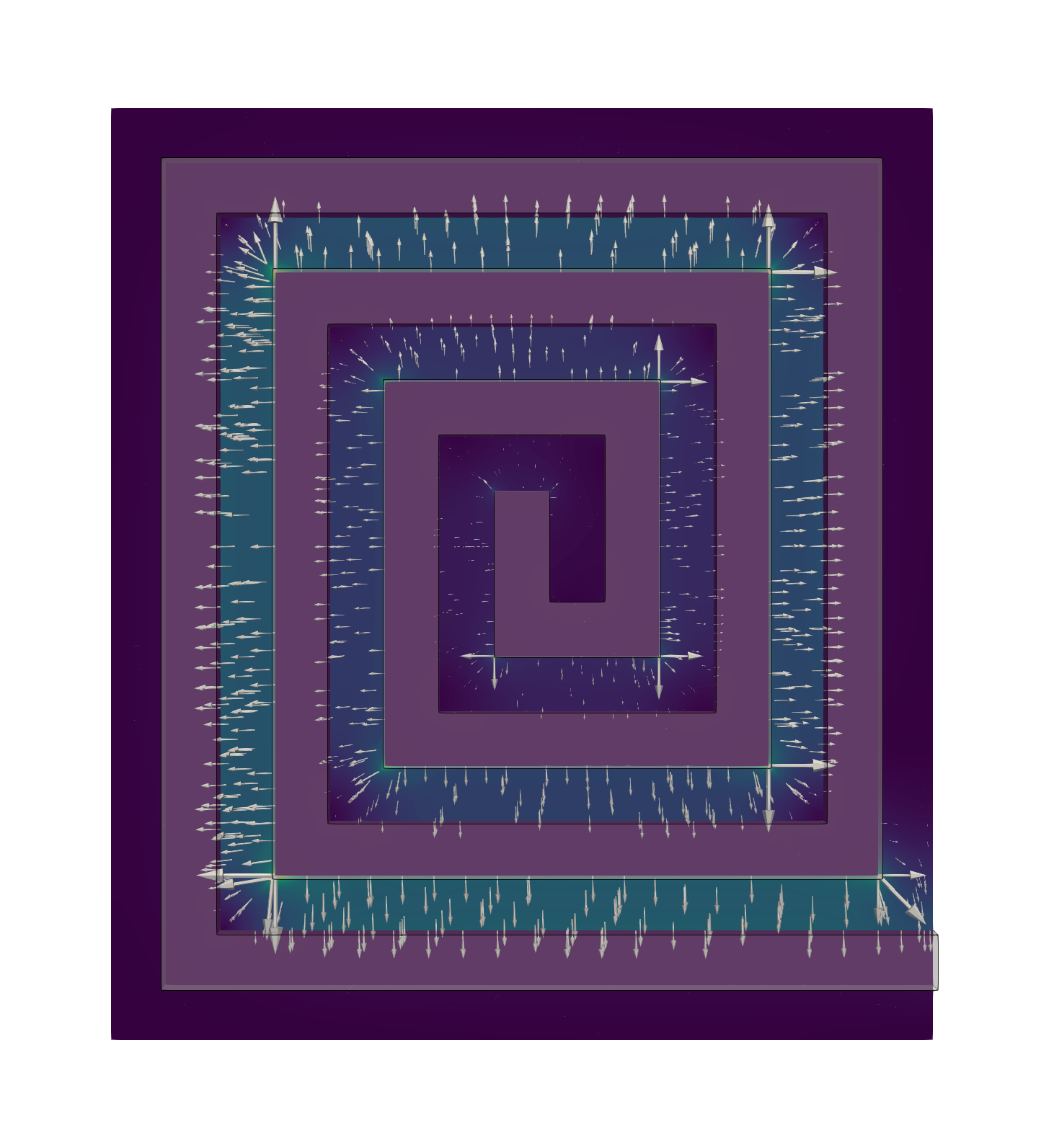}\\
        \hspace{1mm}
        \begin{tikzpicture}[every node/.style={scale=0.8}]
        \pgfplotsset{every tick label/.append style={font=\small}}
            \begin{axis}[%
  hide axis,
  scale only axis,
  colorbar/width=2mm,
  width = 33mm,
  anchor=south,
  point meta min=0.0,
  point meta max=400,
  colormap/viridis,                     %
  colorbar horizontal,                  %
  colorbar sampled,                     %
  colorbar style={
    separate axis lines,
    samples=256,                        %
    xlabel={\small Electric field strength $\big(\SI{}{\volt/\meter}\big)$},
        xtick={0.0,200, 400},
    scaled x ticks=false
  },
]
  \addplot [draw=none] coordinates {(0,0)};
\end{axis}
        \end{tikzpicture}
        \caption{$\norm{\Efield}{2}$, stabilized.}
        \label{fig:Efield_pC}
    \end{subfigure}
    \caption{Top view of the unstabilized and stabilized field solutions at $f = \SI{150}{Hz}$ of the planar coil, sliced in the middle. Note, that for the non-stabilized solutions the magnitudes are scaled logarithmically.}
    \label{fig:planarCoilFields}
\end{figure}

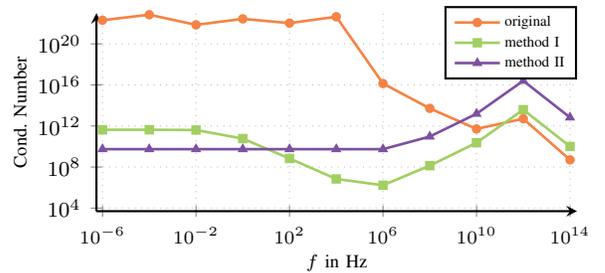
\begin{figure}
    \centering
 \pgfplotstableread[col sep=comma]{data/cond_planarCoil_L3_p2_s2.csv}\cond
                \begin{tikzpicture}
                	\begin{loglogaxis}[tudalineplot, width=0.9\linewidth, height=12em, xlabel={\scriptsize$f$ in Hz}, ylabel={\scriptsize Cond. Number}, xticklabel style={font=\scriptsize}, xmin=5.5e-7,xmax=2e14, ymin=5e3, ymax=1e23, xtick distance=10^(4), ytick distance=10^(4), xlabel shift={-1mm}, ylabel shift={-1mm},legend style={at={(1,0.85)},anchor=east}, legend style={nodes={scale=0.7, transform shape}}]

                 \addplot+[mark size=1.5mm,mark=*,mark options={fill=TUDa-8a,scale=0.3}, TUDa-8a, line width=1pt] table[x index = 0, y index = 1] {\cond};
                \addplot+[mark size=1.5mm,mark=square*,mark options={fill=TUDa-4a,scale=0.3}, TUDa-4a, line width=1pt] table[x index = 0, y index = 3] {\cond};

                    \addlegendentry{original};
                		\addlegendentry{stabilized};
                	\end{loglogaxis}
                \end{tikzpicture}

                    \caption{Condition number of different system matrices for planar coil problem over frequency.}
    \label{fig:condPlanarCoil}
\end{figure}

In \autoref{fig:planarCoilFields} the magnetic flux density and electric field strength are shown. The non-stable original formulation \eqref{eq:discrete1} and \eqref{eq:discrete2} does not yield correct field solutions. Note, that here, due to the large field values, the magnitude is scaled logarithmically. The stabilized versions of the two-step formulation successfully capture not only the induction of the magnetic field but also the capacitance between the coil windings correctly.

\autoref{fig:condPlanarCoil} shows the condition numbers of the different methods over the ordinary frequency. We can again observe that the original system becomes singular for lower frequencies. {The stabilized system {has} a much lower condition number that does not deteriorate for $f \rightarrow \SI{0}{Hz}$. In the static limit $f = \SI{0}{Hz}$, we obtain a condition number of approximately $\SI{1.9e9}{}$.} Again, the original system matrix is singular.

\subsection{Investigation of Inductive Coupling}

\begin{figure}
    \centering
    \begin{subfigure}[B]{0.5\linewidth}
        \centering
        \begin{tikzpicture}
            \node at (0,0) {\includegraphics[width=\linewidth,trim={0 0 1cm 1cm},clip]{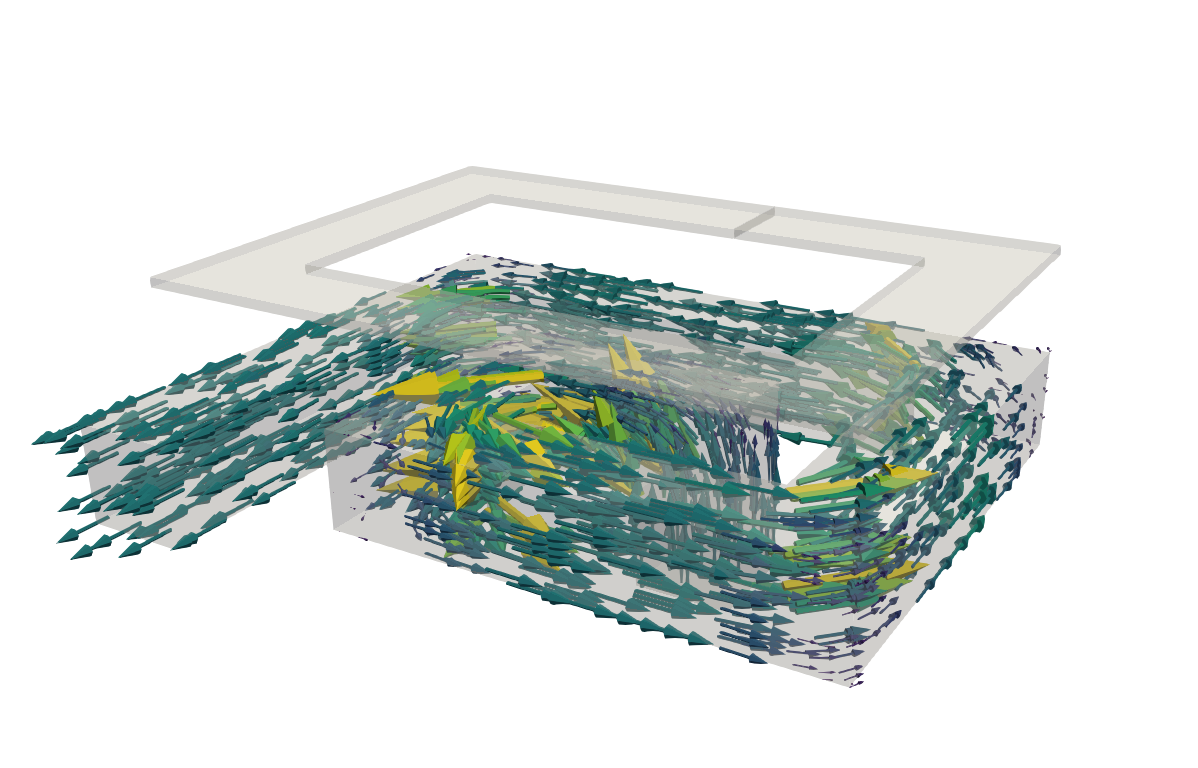}};
            \node at (1,1) {slit};
        \end{tikzpicture}
        \vspace*{-1cm}
        \caption{First step ($\Jfield_\mathrm{e}$).}
        \label{fig:Jstep1}
    \end{subfigure}%
    \begin{subfigure}[B]{0.5\linewidth}
        \centering
        \includegraphics[width=\linewidth,trim={0 0 1cm 1cm},clip]{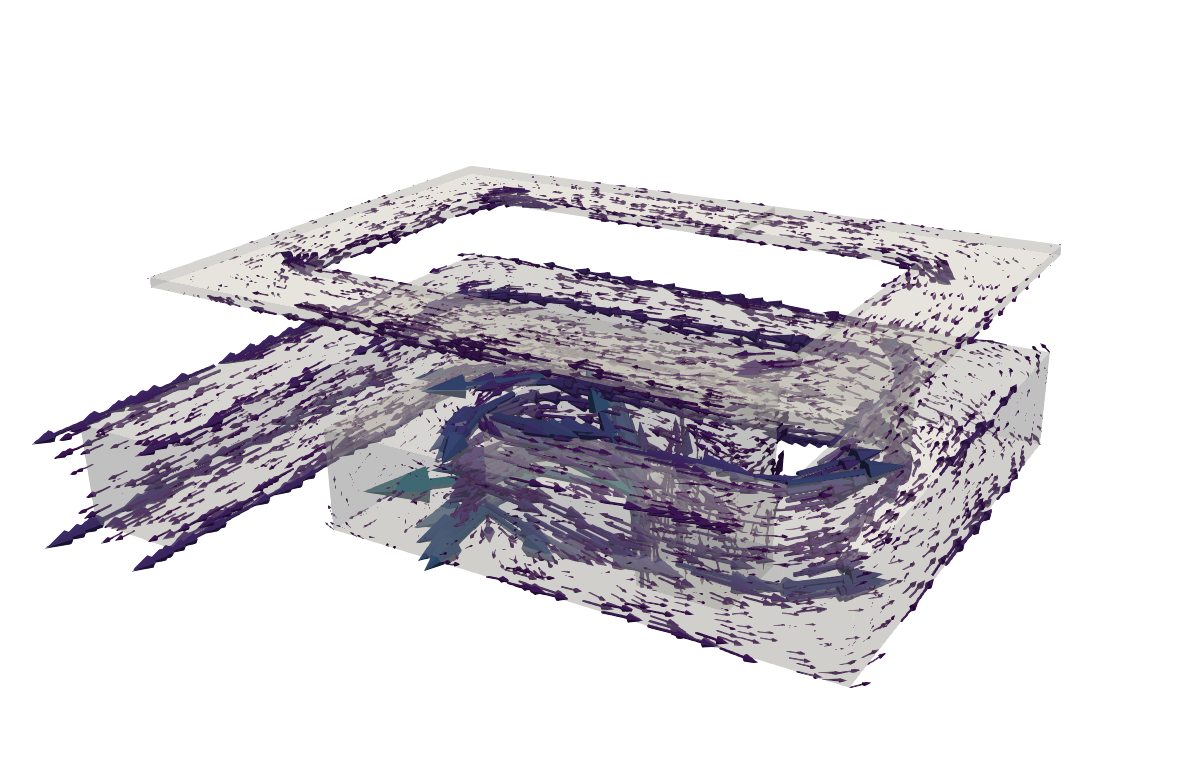}
        \vspace*{-1cm}
        \caption{Both steps ($\Jfield_\mathrm{e} + \Jfield_\mathrm{m}$).}
        \label{fig:Jstep2}
    \end{subfigure}
    \begin{subfigure}[B]{\linewidth}
        \centering
        \begin{tikzpicture}
            \pgfplotsset{every tick label/.append style={font=\small}}
            \begin{axis}[%
                  hide axis,
                  scale only axis,
                  colorbar/width=2mm,
                  width = 50mm,
                  anchor=south,
                  point meta min=0.0,
                  point meta max=1,
                  colormap/viridis,                     %
                  colorbar horizontal,                  %
                  colorbar sampled,                     %
                  colorbar style={
                    separate axis lines,
                    samples=256,                        %
                    xlabel={\small Electric current density $\big(\SI{}{\ampere/\meter^{2}}\big)$},
                    xtick={0.0,0.5, 1},
                    xticklabels={%
                    0.0,
                    $\num{1e8}$,
                    $\num{2e8}$},
                    scaled x ticks=false
                  },
                ]
                \addplot [draw=none] coordinates {(0,0)};
            \end{axis}
        \end{tikzpicture}
    \end{subfigure}
    \caption{Electric current density $\Jfield$ of inductively coupled problem at $f = \SI{150}{Hz}$. Conductor regions are highlighted in gray and the dielectric slit in yellow.}
    \label{fig:trafoExample}
\end{figure}

We investigate the configuration sketched in \autoref{fig:circuit} by utilizing again a planar coil similar to the one from the previous section. The only differences are the number of turns and that we interpret it as the left transformer leg in the context of \autoref{fig:circuit}. To model the right transformer leg, a conducting loop is placed above the coil. We introduce a thin slit in this upper coil with an artificially high relative permittivity of $\permittivity_\mathrm{r} = \SI{7.2e15}{}$ to model the capacitor $C_2$ shown in \autoref{fig:trafoExample}. We employ the same boundary conditions as well as excitation as in \autoref{sec:numPlanarCoil}. {The problem is discretized using 180 patches and {736,054 second order basis functions}.}

\autoref{fig:Jstep1} shows the computed electric current density $\Jfield_\mathrm{e}$ %
after the first step \eqref{eq:discrete2} of the method. Here, we can see that there is no capacitive behavior in the slit of the conductor after step 1 (as expected). Only after completing the second step, shown in \autoref{fig:Jstep2}, the inductive coupling gives rise to a current flow in the loop and the capacitive behavior in the slit allows the electric current  density $\Jfield = \Jfield_\mathrm{e} + \Jfield_\mathrm{m}$
to close through the gap via $i\omega\Dfield_\mathrm{m}$. This underlines the theoretical considerations of \autoref{sec:Darwin} that the capacitive effect cannot be clearly attributed to either step. One can also see that the current amplitude is significantly reduced due to the inductive coupling after step~2.

\section{Conclusion}\label{sec:conclusions}
In this paper, we proposed a novel low-frequency stabilization approach for a two-step formulation of Maxwell’s equations, leveraging a generalized tree-cotree decomposition. Our method ensures numerical stability by maintaining a consistent divergence condition for the magnetic vector potential, preventing a low-frequency breakdown. {A stabilization strategy was introduced and compared, demonstrating its effectiveness in reducing the condition number and preserving accurate field solutions across a broad frequency range including the static cases.

Through numerical experiments, we verified the proposed technique on both academic and application-oriented test cases, confirming its robustness. The stabilized formulation successfully captures inductive and capacitive effects simultaneously.

\section*{Acknowledgment}
The work is supported by the joint DFG/FWF Collaborative Research Centre CREATOR (DFG: Project-ID 492661287/TRR 361; FWF: 10.55776/F90) at TU Darmstadt, TU Graz and JKU Linz. We acknowledge the funding of The "Ernst Ludwig Mobility Grant" of the Association of Friends of Technical University of Darmstadt e.V.

\vspace{2em}


\begin{thebibliography}{10}
\providecommand{\url}[1]{#1}
\csname url@samestyle\endcsname
\providecommand{\newblock}{\relax}
\providecommand{\bibinfo}[2]{#2}
\providecommand{\BIBentrySTDinterwordspacing}{\spaceskip=0pt\relax}
\providecommand{\BIBentryALTinterwordstretchfactor}{4}
\providecommand{\BIBentryALTinterwordspacing}{\spaceskip=\fontdimen2\font plus
\BIBentryALTinterwordstretchfactor\fontdimen3\font minus
  \fontdimen4\font\relax}
\providecommand{\BIBforeignlanguage}[2]{{%
\expandafter\ifx\csname l@#1\endcsname\relax
\typeout{** WARNING: IEEEtran.bst: No hyphenation pattern has been}%
\typeout{** loaded for the language `#1'. Using the pattern for}%
\typeout{** the default language instead.}%
\else
\language=\csname l@#1\endcsname
\fi
#2}}
\providecommand{\BIBdecl}{\relax}
\BIBdecl

\bibitem{Monk_2003aa}
P.~Monk, \emph{\BIBforeignlanguage{english}{Finite Element Methods for
  {Maxwell}'s Equations}}.\hskip 1em plus 0.5em minus 0.4em\relax Oxford:
  Oxford University Press, 2003.

\bibitem{Zhu_2010aa}
J.~Zhu and D.~Jiao, ``A theoretically rigorous full-wave finite-element-based
  solution of {Maxwell}’s equations from dc to high frequencies,''
  \emph{{IEEE} Trans. Adv. Packag.}, vol.~33, no.~4, pp. 1043--1050, 11 2010.

\bibitem{Albanese_1990aa}
R.~Albanese and G.~Rubinacci, ``\BIBforeignlanguage{english}{Magnetostatic
  field computations in terms of two-component vector potentials},''
  \emph{\BIBforeignlanguage{english}{Int. J. Numer. Meth. Eng.}}, vol.~29,
  no.~3, pp. 515--532, 1990.

\bibitem{Munteanu_2002aa}
\BIBentryALTinterwordspacing
I.~Munteanu, ``\BIBforeignlanguage{english}{Tree-cotree condensation
  properties},'' \emph{\BIBforeignlanguage{english}{{ICS} Newslett.}}, vol.~9,
  pp. 10--14, 2002. [Online]. Available:
  \url{https://www.compumag.org/wp/newsletter/}
\BIBentrySTDinterwordspacing

\bibitem{Dyczij-Edlinger_1999aa}
R.~Dyczij-Edlinger, G.~Peng, and J.-F. Lee, ``Efficient finite element solvers
  for the {Maxwell} equations in the frequency domain,'' \emph{Comput. Meth.
  Appl. Mech. Eng.}, vol. 169, no. 3–4, pp. 297--309, 02 1999.

\bibitem{Hiptmair_2008aa}
R.~Hiptmair, F.~Kramer, and J.~Ostrowski, ``A robust {Maxwell} formulation for
  all frequencies,'' \emph{{IEEE} Trans. Magn.}, vol.~44, no.~6, pp. 682--685,
  06 2008.

\bibitem{Jochum_2015aa}
M.~T. Jochum, O.~Farle, and R.~Dyczij-Edlinger, ``A new low-frequency stable
  potential formulation for the finite-element simulation of electromagnetic
  fields,'' \emph{{IEEE} Trans. Magn.}, vol.~51, no.~3, p. 7402304, 03 2015.

\bibitem{Jochum_2016aa}
M.~Jochum, O.~Farle, and R.~Dyczij-Edlinger, ``\BIBforeignlanguage{english}{A
  symmetric and low-frequency stable potential formulation for the
  finite-element simulation of electromagnetic fields},'' in
  \emph{\BIBforeignlanguage{english}{Scientific Computing in Electrical
  Engineering {SCEE} 2014}}, ser. Mathematics in Industry, no.~23.\hskip 1em
  plus 0.5em minus 0.4em\relax Berlin: Springer, 05 2016, pp. 63--71.

\bibitem{Eller_2017aa}
M.~Eller, S.~Reitzinger, S.~Schöps, and S.~Zaglmayr,
  ``\BIBforeignlanguage{english}{A symmetric low-frequency stable broadband
  {Maxwell} formulation for industrial applications},''
  \emph{\BIBforeignlanguage{english}{{SIAM} J. Sci. Comput.}}, vol.~39, no.~4,
  pp. B703--B731, 08 2017.

\bibitem{Ho_2016aa}
S.~L. Ho, Y.~Zhao, W.~N. Fu, and P.~Zhou, ``Application of edge elements to 3-d
  electromagnetic field analysis accounting for both inductive and capacitive
  effects,'' \emph{IEEE Transactions on Magnetics}, vol.~52, no.~3, pp. 1--4,
  2016.

\bibitem{Ostrowski_2021aa}
J.~Ostrowski and R.~Hiptmair, ``Frequency-stable full {Maxwell} in
  electro-quasistatic gauge,'' \emph{{SIAM} J. Sci. Comput.}, vol.~43, no.~4,
  pp. B1008--B1028, 01 2021.

\bibitem{Formaggia_2012aa}
L.~Formaggia, F.~Saleri, and A.~Veneziani, \emph{Solving Numerical {PDEs}:
  Problems, Applications, Exercises}.\hskip 1em plus 0.5em minus 0.4em\relax
  Springer, 2012.

\bibitem{Clemens_2022aa}
\BIBentryALTinterwordspacing
M.~Clemens, M.-L. Henkel, F.~Kasolis, M.~Günther, H.~De~Gersem, and
  S.~Schöps, ``\BIBforeignlanguage{english}{Electromagnetic quasistatic field
  formulations of {Darwin} type},'' \emph{\BIBforeignlanguage{english}{{ICS}
  Newslett.}}, vol.~29, no.~3, pp. 3--9, 03 2022, arxiv:2204.06286. [Online].
  Available: \url{https://www.compumag.org/wp/newsletter/}
\BIBentrySTDinterwordspacing

\bibitem{Balian_2023aa}
D.~Balian, M.~Merkel, J.~Ostrowski, H.~De~Gersem, and S.~Schöps,
  ``Low-frequency stabilization of dielectric simulation problems with
  conductors and insulators,'' \emph{{IEEE} Trans. Dielectr. Electr. Insul.},
  vol.~30, no.~6, 12 2023, arxiv:2302.00313.

\bibitem{Kaimori_2024aa}
H.~Kaimori, T.~Mifune, A.~Kameari, and S.~Wakao, ``Low-frequency stabilized
  formulations of {Darwin} model in time-domain electromagnetic finite-element
  method,'' \emph{{IEEE} Trans. Magn.}, vol.~60, no.~3, pp. 1--5, 2024.

\bibitem{Buffa_2011aa}
A.~Buffa, J.~Rivas, G.~Sangalli, and R.~Vázquez~Hernández,
  ``\BIBforeignlanguage{english}{Isogeometric discrete differential forms in
  three dimensions},'' \emph{\BIBforeignlanguage{english}{{SIAM} J. Numer.
  Anal.}}, vol.~49, no.~2, pp. 818--844, 2011.

\bibitem{Albanese_1988aa}
R.~Albanese and G.~Rubinacci, ``\BIBforeignlanguage{english}{Integral
  formulation for {3D} eddy-current computation using edge elements},''
  \emph{\BIBforeignlanguage{english}{{IEE} Proc. Sci. Meas. Tech.}}, vol. 135,
  no.~7, pp. 457--462, 09 1988.

\bibitem{Manges_1995aa}
J.~B. Manges and Z.~J. Cendes, ``A generalized tree-cotree gauge for magnetic
  field computation,'' \emph{{IEEE} Trans. Magn.}, vol.~31, no.~3, pp.
  1342--1347, 1995.

\bibitem{Rapetti_2022aa}
F.~Rapetti, A.~Alonso~Rodríguez, and E.~De~Los~Santos, ``On the tree gauge in
  magnetostatics,'' \emph{J}, vol.~5, no.~1, pp. 52--63, 2022.

\bibitem{Vazquez_2016aa}
R.~Vázquez, ``A new design for the implementation of isogeometric analysis in
  {Octave} and {Matlab}: {GeoPDEs} 3.0,'' \emph{Comput. Math. Appl.}, vol.~72,
  no.~3, pp. 523--554, 08 2016.

\bibitem{Herles_2025ab}
\BIBentryALTinterwordspacing
L.~Herles, ``{LowFrequencyStableFullMaxwell},'' 02 2025. [Online]. Available:
  \url{https://www.doi.org/10.5281/zenodo.14810885}
\BIBentrySTDinterwordspacing

\bibitem{Hiptmair_2000aa}
R.~Hiptmair, ``\BIBforeignlanguage{english}{Multilevel gauging for edge
  elements},'' \emph{\BIBforeignlanguage{english}{Computing}}, vol.~64, no.~2,
  pp. 97--122, 2000.

\end{thebibliography}
\end{document}